\documentclass[review]{elsarticle}

\usepackage{lineno}
\modulolinenumbers[1]

\journal{Journal of \LaTeX\ Templates}

\usepackage{amsmath,color,colordvi,caption,subcaption,url}
\usepackage{graphicx}
\usepackage{amssymb}
\usepackage{multirow}
\usepackage{picinpar}
\usepackage{algorithm}
\usepackage{algorithmic}
\usepackage{footmisc}
\usepackage{array}
 
 \usepackage[colorlinks,
            linkcolor=red,
            anchorcolor=green,
            citecolor=blue
            ]{hyperref}
  \usepackage{cleveref}          
\graphicspath{ {./images/} }

\def\bd{\mathbf{d}}

\def\bw{\mathbf{w}}
\def\bx{\mathbf{x}}
\def\by{\mathbf{y}}
\def\bz{\mathbf{z}}
\def\bg{{\boldsymbol g}}

\def\supp{\mathrm{supp}}

\def\time{\texttt{Time}}
\def\ser{\texttt{SER}}

\def\rx{\mathrm{x}}
\def\NSLR{\texttt{NSLR}}
\def\GPGN{\texttt{GPGN}}
\def\NTGP{\texttt{NTGP}}
\def\GraSP{\texttt{GraSP}}
\def\LARS{\texttt{LARS}}
\def\GIST{\texttt{GIST}}
\def\APG{\texttt{APG}}
\def\SLEP{\texttt{SLEP}}

\newtheorem{proposition}{Proposition}[section]
\newtheorem{theorem}[proposition]{Theorem}
\newtheorem{lemma}[proposition]{Lemma}

\newtheorem{definition}[proposition]{Definition}
\newtheorem{remark}[proposition]{Remark}
\newtheorem{example}[proposition]{Example}

\numberwithin{equation}{section}
\begin{document}

\begin{frontmatter}

\title{ {An Extended Newton-type Algorithm for $\ell_2$-Regularized Sparse Logistic Regression  and Its Efficiency
for Classifying Large-scale Datasets}}


\author[mymainaddress]{Rui Wang}
\ead{wangruibjtu@bjtu.edu.cn}

\author[mymainaddress]{
        Naihua Xiu}
\ead{nhxiu@bjtu.edu.cn}

\author[mysecondaryaddress]{Shenglong Zhou\corref{mycorrespondingauthor}}
\cortext[mycorrespondingauthor]{Corresponding author}
\ead{slzhou2021@163.com}

\address[mymainaddress]{Department of Applied Mathematics, Beijing Jiaotong University, Beijing, China}
\address[mysecondaryaddress]{Department of EEE, Imperial College London, London, UK}

\begin{abstract}
Sparse logistic regression, {as an effective tool of classification,} has been developed tremendously in recent two decades, from its origination the $\ell_1$-regularized version to the sparsity constrained models. This paper is carried out on the sparsity constrained logistic regression by the Newton method. We begin with establishing its first-order optimality condition associated with a $\tau$-stationary point. This point can be equivalently interpreted as a system of equations which is then efficiently solved by the Newton method. The  method has a considerably low computational complexity and enjoys global and quadratic convergence properties. Numerical experiments on random and real data demonstrate its superior performance when against  seven state-of-the-art solvers.
\end{abstract}


\end{frontmatter}

\section{Introduction}\label{sec:introduction}
 As one of effective tools of classification,  logistic regression has its high reputation with extensive applications ranging from machine learning, data mining,  pattern recognition,  medical science to statistics. It describes the relationship between a sample data  $\bx $  and its associated binary response/label $y \in \{0,1\}$  through the  conditional probability
\begin{eqnarray}\label{LR-cond-prob}{\rm Pr}\left(y \left|\right. \bx, \bz\right) =  \left({1+e^{-\langle \bx ,\bz \rangle}} \right)^{-1}, \end{eqnarray}
where ${\rm Pr}(y  | \bx,\bz)$ is the conditional probability of the label $y$, given the sample $\bx$ and a parameter vector $\bz$, and $\langle \bx ,\bz \rangle$ is the vector inner product. To find the maximum likelihood estimate of the parameter $\bz$,  a set of  {$n$ i.i.d. (independently and identically distributed)} samples $\{(\bx_i,y_i), i=1,2,\ldots,n\}$ are first drawn, where  $\bx_i \in \mathbb{R}^{p}$ and $y_i \in \{0,1\}$, yielding a joint likelihood of the interested parameter/classifier $\bz$.  Then the maximum likelihood estimate is obtained by minimizing the classical logistic regression loss function,
\begin{eqnarray}\label{LR}
 \ell(\bz):=\frac{1}{n}\sum_{i=1}^{n}\left(\ln\left(1+ e^{\langle\bx_i,\bz\rangle}\right)-y_i\langle\bx_i,\bz\rangle\right).
\end{eqnarray}
The logistic loss function is strictly convex and thus admits a unique minimizer provided that the sample matrix is full row rank. Therefore, the minimization performs relatively well when the number of samples is larger than the number of features, i.e.,  $n\geq p$.  But, the case $n<p$ may lead to an over-fitting: the solved classifier through minimizing (\ref{LR}) well fits the model (making the loss sufficiently small) on training data but behaves poorly on unseen data.

On the one hand,  the case  $n<p $ occurs  often in many real applications. For instance, one piece of gene expression data sample is made of thousands of genes whilst common medical equipments are only able to obtain very limited samples. In image processing, an image consists of large amounts of pixels, which is far more than the number of observed images. One the other hand, despite numerous features in those data, there is only a small portion that is of importance. For example, apart from the classification task, the micro-array data experiments also attempt to identify a small set of informative genes (to distinguish the tumour and the normal tissues) in each gene expression data so as to remove the irrelevant genes to simplify the inference. This naturally gives rise to the topic of the sparse logistic regression.
\subsection{Sparse logistic regression}
Sparse logistic regression (SLR) was originated from the $\ell_1$-regularized logistic regression in \cite{tibshirani1996regression},
\begin{equation}\label{SLR-L1}
\min_{\bz\in \mathbb{R}^p }~~ \ell(\bz)+\nu\|\bz\|_1,
\end{equation}
where $\|\bz\|_1$ is the $\ell_1$-norm and $\nu>0$. Under the help of $\ell_1$-regularization, this model is capable of rendering a sparse solution allowing for capturing key features among others. A vector is called sparse if only a few entries are non-zero  and the rest are zeros.  With the advance in sparse optimization in recent decade, (\ref{SLR-L1}) has been extensively extended to the following general model,
\begin{equation}\label{SLR-p}
\min_{\bz\in \mathbb{R}^p }~~\ell_\phi(\bz):=\ell(\bz)+\phi_{\nu}(\bz),
\end{equation}
where the regularized function $\phi_{\nu}(\bz):\mathbb{R}^p\rightarrow\mathbb{R}$ is designed to pursue a sparse solution and associated with some given non-negative parameters $\nu$.

An alternative is to consider logistic regression with a sparsity constraint, which was first studied in \cite{bahmani2013greedy,plan2013robust} separately and then well investigated in \cite{Wang}. They perform the following sparsity constrained logistic regression
\begin{equation}\label{SLR}
\min_{\bz\in \mathbb{R}^p }~~ \ell(\bz) ,~~{\rm s.t.}~~ \|\bz\|_0\leq s,
\end{equation}
where $\|\bz\|_0$ is the $\ell_0$ pseudo norm of $\bz$, counting the number of non-zero elements of $\bz$. The discreteness of the sparsity constraint makes tackling this model NP-hard. Nevertheless, compared with the regularized model, the sparsity constrained version enjoys  various appealing features, such as being penalty parameter-free, ease of sparsity controlling, and low computational complexity in terms of numerical computation and so  forth.

Therefore, the generalization of the problem \eqref{SLR}, where $\ell(\bz) $ is replaced by a more general function,  has been thoroughly investigated in \cite{beck15,pan2015solutions} since it was first introduced by \cite{bahmani2013greedy} and \cite{Beck13}. Particularly, in statistics, the model with the logistic loss function being replaced by the least squares of linear regression is the so-called best subspace/feature selection \cite{hastie2017extended, mazumder2017subset, hazimeh2018fast, xie2018ccp, Pang2019}. Those research bring fruitful results and provide  a series of effective numerical tools to conquer the NP-hardness.  

However, as stated in \cite{bahmani2013greedy} that `\textit{one can achieve arbitrarily small loss values by tending the parameters to infinity along certain directions}' for  \eqref{SLR}, authors \cite{bahmani2013greedy} suggests to address the following regularized model
\begin{equation}\label{SLR-L2}
\min~~f(\bz):=\ell(\bz)+(\lambda/2)\|\bz\|_2^2,~~{\rm s.t.}~~\|\bz\|_0\leq s,
\end{equation}
where $\lambda>0$ is a given penalty parameter. Now the objective function $f$ is strongly convex and thus \eqref{SLR-L2} admits {finitely} many (local or global) bounded minimizers. So the work in this paper is carried out along with this model.

\subsection{Methods of solving SLR}\label{methods-SLR}
Since there is a vast body of methods that have been proposed to deal with the sparse optimization problems containing the  SLR as a special case, we present a brief overview of methods that process the problems  (\ref{SLR-p})-(\ref{SLR-L2}) directly. 

\paragraph{Regularization methods} Most versions of the model (\ref{SLR-p}) are unconstrained and continuous. Then generic optimization methods, known as the relaxation (regularization) methods from the perspective of optimization, are tractable. Dependent on the convexity of the penalty functions $\phi _{\nu}$,  those methods  can be summarized into two {categories}.

 {Convex regularizations} are mainly associated with the usage of $\ell_1$-norm: 
\begin{itemize}
\item $\phi _{\nu}(\bz)=\nu\|\bz\|_1$. Some earliest work can be traced back to    \cite{figueiredo2003adaptive, krishnapuram2005sparse}, where expectation maximization methods were developed. Later relevant work can be found in \cite{andrew2007scalable, koh2007interior,yu2010quasi,shi2010fast,yuan2010comparison}.
  
\item  $\phi_{\boldsymbol\nu}(\bz)=\nu_1\|\bz\|^2_2+\nu_2\|\bz\|_1$, where $\boldsymbol \nu:=(\nu_1, \nu_2)>0$. For this penalty, two powerful packages  \texttt{SLEP} \cite{liu2009slep} and  \texttt{GLMNET} \cite{friedman2010regularization, yuan2012improved} have been created. 

\item $\phi _{\nu}(\bz)=\nu\|\bz\|^2+\delta_{\|\bz\|_1\leq t}(\bz)$, where $t>0$  {is} a given parameter, and $\delta_{\|\bz\|_1\leq t}(\bz)=0$ if $\|\bz\|_1\leq t$ and $+\infty$ otherwise.  
Such a problem can be addressed by  \texttt{Lassplore}  \cite{liu2009large} or  \texttt{SLEP} \cite{liu2009slep}.
 When $\nu=0$, the above model is the $\ell_1$ constrained logistic regression, which was addressed by \texttt{IRLS-LARS} in \cite{lee2006efficient}. Here, \texttt{LARS} was adopt from \cite{efron2004least}.
 \end{itemize}
 
 {Nonconvex regularizations} differ slightly. In the early stage, scholars from statistics have proposed a number of excellent methods including the smoothly clipped absolute deviation (\texttt{SCAD}  \cite{fan2001variable}), one step local linear approximation \cite{zou2008one} and the group bridge method for multiple regression problems \cite{huang2009group}. Then, a  general iterative shrinkage and thresholding algorithm (\texttt{GIST}) has been proposed in \cite{gong2013general}. 
Recently, the accelerated proximal gradient method (\texttt{APG}) in  \cite{li2015accelerated}, the efficient hybrid optimization algorithm
for non-convex regularized problems (\texttt{HONOR}) in \cite{gong2015honor}  and the proximal Newton method based on the scheme of the difference of two convex functions in
\cite{rakotomamonjy2016dc} are worth exploring.
 
\paragraph {Greedy Methods} An impressive body of  approaches have been developed to solve the sparsity constrained models (\ref{SLR}) or \eqref{SLR-L2}. The first work in \cite{bahmani2013greedy} generalized the compressive sampling matching pursuit \cite{needell2009cosamp} to derive the gradient support pursuit (\texttt{GraSP}).  Then authors in \cite{lozano2011group} adopted the orthogonal matching pursuit (\texttt{OMP} \cite{mallat1993matching})  to develop a group \texttt{OMP} method. Other relating methods can be seen those in \cite{lu2013sparse, yuan2014gradient,pan2017convergent}. Very lately, three effective Newton type methods have been designed. They are the Newton greedy pursuit method \texttt{NTGP} in \cite{yuan2017newton}, greedy projected gradient-Newton  method (\texttt{GPGN} \cite{Wang}) and the fast Newton hard thresholding pursuit  \cite{chen2017fast}. {In particular, we would like to mention the  methods, the zero-CW search method and the full-CW search method, proposed in \cite{beck15}. Both methods first carefully search an index set $T$ and then solve a subproblem where the variable has support within $T$ to update the next point. }

\subsection{Our contributions}
Those aforementioned methods have been testified to have the excellent numerical performance to deal with (\ref{SLR}) or (\ref{SLR-L2}). However, only a very few of them established strong theoretical guarantees (such as global convergence property or quadratic convergence rate) from the perspective of deterministic optimization. Therefore, in this paper, we aim to develop a second-order method that possesses such strong theoretical guarantees.  The main contributions are summarized as follows.

{\rm C1)}  We start with establishing the optimality condition of the model (\ref{SLR-L2}) by introducing a  $\tau$-stationary point (see Definition \ref{alpha} for more details) which turns out to be at least a locally optimal solution by Theorem \ref{relation1}. More importantly, a $\tau$-stationary point draws forth a system of equations (\ref{equation}) that makes the classic Newton method applicable.

{\rm C2)} Differing with any of the above mentioned algorithms, we perform the Newton method on solving a system of equations (\ref{equation}), one of the optimality conditions of the problem (\ref{SLR-L2}).  {The proposed Newton method for SLR (\NSLR\ for short) has a simple framework} (see    \Cref{Alg-NSLR}) that makes its implementation easy and a low computational complexity per each iteration. Such a low computational complexity is due to a small-scale linear equation system with  $s$ variables and $s$ equations being solved to update the Newton direction.

{\rm C3)}  {It is worth mentioning that the standard Newton-type methods derive the directions for a fixed system of equations. However,  in each iteration, the system of equations \eqref{equation} varies when the index set $\alpha$ changes. Consequently, \NSLR\  updates Newton directions on unfixed systems  of equations. Because of this, some common approaches to establish the convergence results of Newton-type methods for solving a fixed system of equations fail to be employed for \NSLR.} Nevertheless, we still show that the whole sequence generated by \NSLR\ converges to a $\tau$-stationary point, at least a locally optimal solution. Moreover, the convergence enjoys a quadratic rate, well testifying the proposed method would perform extraordinarily theoretically. 
 
{\rm C4)} Finally, the efficiency of \NSLR\ is demonstrated against seven state-of-the-art methods by solving a number of randomly generated and real datasets. The fitting accuracy and computational speed are very competitive. Especially, in high dimensional data setting,  \NSLR\ outperforms the others in terms of the computational time.

 {We note that there are some methods that also have a close link to $\tau$-stationary point, such as  two methods in \cite{beck15} and \texttt{GPGN} \cite{Wang}. We would like to highlight the difference between them and \NSLR.  For the methods in \cite{beck15},  since an optimal solution to a subproblem needs to be found to update the next point in each step,  the methods can terminate within finitely many steps. However, \NSLR\ updates the next point by Newton direction with a line search scheme and has been shown to enjoy the global and quadratic convergence properties. For \texttt{GPGN}, the procedure {\tt IHT} from \cite{Beck13} to update the next point dominates most steps, and Newton steps are imposed only when two consecutive points have the same support sets. It is shown to converge  quadratically only when the solution has $s$ nonzeros but sublinearly otherwise.  By contrast, \NSLR\ always performs Newton step to update the next point and converges quadratically without additional assumptions.   Moreover, differing from papers \cite{beck15} and \cite{Wang} where comprehensive optimality conditions have been investigated, the primary aim of this paper is to develop a Newton-type method and establish its convergence properties.}
\subsection{Organization and notation}
This paper is organized as follows. To explore the optimality conditions of (\ref{SLR-L2}), the next section introduces the $\tau$-stationary point by Definition \ref{alpha} and establishes  its relationships with a local/global minimizer in Theorem \ref{relation1}. This $\tau$-stationary point is then  equivalently transferred to a system of equations  (\ref{equation}). Section \ref{Section3} develops the method
\NSLR, an abbreviation for Newton method for SLR, which turns out to have a simple algorithmic framework and low computational complexity.  The global and quadratic convergence properties of the method are then established. In Section \ref{Section5}, the superior performance of \NSLR\ is demonstrated against some of the state-of-the-art solvers on randomly generated and real datasets in high dimensional scenarios. Concluding remarks are made in the last section.

 We end this section by defining some notation employed throughout this paper. Let $X:=(\bx_1,\bx_2,\ldots,\bx_n)^T\in \mathbb{R}^{n\times p  }$ be the sample matrix and $\by=(y_1,y_2,\ldots,y_n)^T \in \mathbb{R}^{n}$ be the response vector.
 For an index set $\alpha\subseteq  [p]:=\{ 1,2,\ldots,p\}$, let $|\alpha|$ be the cardinality of $\alpha$ and $\overline{\alpha}:=[p]\setminus \alpha$ be the complementary set of $\alpha$. The support  set of a vector $\bz$ is denoted by $\supp(\bz):=\{i\in[p]: z_i\neq0\}$. We denote $[\bz]^{\downarrow}_i$ the $i$th largest (in absolute) elements of $\bz$. Write $\bz_\alpha\in\mathbb{R}^{|\alpha|}$   as the sub-vector  of $\bz$ containing elements indexed on $\alpha$. Similarly, for a matrix $A\in\mathbb R^{p\times p}$, $A_{\alpha\beta}$ is the sub-matrix containing rows indexed on ${\alpha}$ and columns indexed on $\beta$, particularly, $A_{\alpha:}=A_{\alpha[p]}$ if $\beta=[p]$.  
Let $\|\cdot\|$  denote the Spectral norm for a matrix and Euclidean norm for a vector respectively. Furthermore, $\lambda_{\min}(A)$ and $\lambda_{\max}(A)$ are the minimal and maximal eigenvalues of  $A$. 

 \section{Optimality} \label{Section2}

This section is devoted to investigate the optimality conditions of (\ref{SLR-L2}), before which we summarize some properties of   $\ell(\bz)$  from \cite{Wang}.
\begin{proposition}[Lemma 2.2-2.4, Lemma A.3 \cite{Wang}]\label{basic-pro}
The function $\ell(\bz)$ is twice continuously differentiable  and  has the following basic properties:
\begin{itemize}
\item[i)] It is non-negative, convex and strongly smooth on $\mathbb{R}^p$ with a parameter $$\lambda_{\rx}:= \lambda_{\max}(X^T X)/(4n).$$ 
\item[ii)]The gradient is Lipschitz continuous with the Lipschitz constant $\lambda_{\rx}$.
\item[iii)] The Hessian matrix  is Lipschitz continuous with the Lipschitz constant $M:= 12\lambda_{\rx}{\max}_{i\in[n]}\|\bx_i\|_1$.
\end{itemize}
\end{proposition}
These properties of $\ell(\bz)$ are also enjoyed by the function $f(\bz)=\ell(\bz)+(\lambda/2)\|\bz\|^2$. Since the proofs are easy, we only summarize them here. The function $f$ is strongly convex
with a constant $\lambda$, and strongly smooth with a parameter $L:=\lambda+\lambda_{\rx}$, namely, for any $\bz,\bz'\in\mathbb{R}^p$,
  \begin{eqnarray}\label{str-con}
  \begin{array}{lll}
  f(\bz) &\geq&  f(\bz')+\langle \nabla f(\bz'),\bz-\bz' \rangle + ({\lambda }/{2})   \|\bz-\bz'\|^2,\\
f(\bz) &\leq&  f(\bz')+\langle \nabla f(\bz'),\bz-\bz' \rangle + ({L}/{2})   \|\bz-\bz'\|^2.
  \end{array}
\end{eqnarray}
The gradient is Lipschitz continuous with the Lipschitz constant $L$, namely,
  \begin{eqnarray*}
\| \nabla f(\bz)-\nabla f(\bz')\| \leq L  \|\bz-\bz'\|^2, 
\end{eqnarray*}
 for any $\bz,\bz'\in\mathbb{R}^p$. The Hessian matrix is Lipschitz continuous with the Lipschitz constant $M$, namely, for any $\bz,\bz'\in\mathbb{R}^p$,
  \begin{eqnarray}\label{hess-lip}
\| \nabla^2 f(\bz)-\nabla^2 f(\bz')\| \leq M \|\bz-\bz'\|^2.
\end{eqnarray}
When it comes to characterize the solutions of the problem $(\ref{SLR-L2})$, we need the projection of a vector $\bz \in \mathbb{R}^p$ onto the feasible region defined by
$$\Pi_s(\bz):={\rm argmin}_{\bx \in \mathbb{R}^p}\{\|\bz-\bx\|: \|\bx\|_0\leq s\},$$
which sets all but $s$ largest absolute value components of $\bz$ to zero. Since  the right hand side may have multiple solutions,  $\Pi_s(\bz)$ is a set. Based on the projection, we introduce the concept of the $\tau$-stationary point {which is also known as the $L$ stationary point in \cite[Definition 2.3]{Beck13}}.  
\begin{definition}\label{alpha}{\cite[Definition 2.3]{Beck13}}
A point $\bz$ is called a $\tau$-stationary point of the problem  $\eqref{SLR-L2}$
if there is a $\tau> 0$ satisfying
\begin{equation}\label{alphaequation}
 \bz \in \Pi_s\left(\bz-\tau \nabla f(\bz)\right).
\end{equation}
\end{definition}
By \cite[Lemma 2.2 ]{Beck13}, $\bz$ is a $\tau$-stationary point if and only if  
\begin{eqnarray}\label{agradient}
\|\bz\|_{0} \leq s,~~~~\tau\left|(\nabla f(\bz))_i\right| \left\{
             \begin{array}{lll}
             = & 0, &i \in \supp(\bz),\\
             \leq & [\bz]^{\downarrow}_s, &i \notin \supp(\bz).
             \end{array}
        \right.
\end{eqnarray}
Based on the definition of the $\tau$-stationary point, our first main result is establishing its relationships with a locally/globally optimal solution to   (\ref{SLR-L2}).
 \begin{theorem}\label{relation1}
 The following results hold for the problem $(\ref{SLR-L2})$.
\begin{itemize}
    \item[{\rm i)}] A global minimizer is a $\tau$-stationary point $\bz^*$ for any $ 0<\tau < 1/L.$
    \item[{\rm ii)}] A $\tau$-stationary point $\bz^*$ for some $\tau>0$ is a unique local minimizer if $\|\bz^*\|_0= s$ and a unique global minimizer if $\|\bz^*\|_0 < s$.
    \item[{\rm iii)}] A  $\tau$-stationary point for some $\tau>1/\lambda$  is a unique global minimizer.
\end{itemize}
\end{theorem}
{\bf Proof}~i)~ The proof is the same as that in \cite[Theorem 2.2]{Beck13}.

ii)~ Let $\bz^*$ be a  $\tau$-stationary point for some $\tau >0$. Then we have (\ref{agradient}), i.e., \begin{eqnarray}\label{agradient-*}
\tau\left|(\nabla f(\bz^*))_i\right| \left\{
             \begin{array}{lll}
             = & 0, &i \in \supp(\bz^*)=:\alpha_*,\\
             \leq & [\bz^*]^{\downarrow}_s, &i \notin \supp(\bz^*).
             \end{array}
        \right.
\end{eqnarray}
For the case $\|\bz^*\|_0=s$, consider a local region $N(\bz^*):=\{\bz: \|\bz-\bz^*\|<[\bz^*]^{\downarrow}_s\}$. Then for any feasible point $\bz \in N(\bz^*)$ and any $i\in\alpha_*$, we have $|z_i|\geq |z^*_i|-|z^*_i-z_i|> |z^*_i| - [\bz^*]^{\downarrow}_s \geq0,$
 which means $\alpha_*\subseteq \supp(\bz)$. Since $\|\bz\|_0\leq s=|\alpha_*|$, it holds $
\alpha_* = \supp(\bz)$
for any feasible point $\bz \in N(\bz^*)$, namely $\bz_{\overline{\alpha}_*}=\bz^*_{\overline{\alpha}_*}=0$.
Then the strong convexity of $f$ in \eqref{str-con} leads to
 \begin{eqnarray}\label{stongly-con-1} \arraycolsep=1.4pt\def\arraystretch{1.25}
 \begin{array}{lll}    
2f(\bz)-2f(\bz^*) 
&\geq &2\langle \nabla f(\bz^*),\bz-\bz^* \rangle +  {\lambda} \|\bz-\bz^*\|^2\\
&=& 2 \left\langle (\nabla f(\bz^*))_{\overline{\alpha}_*},(\bz-\bz^*)_{\overline{\alpha}_*} \right\rangle  + {\lambda} \|\bz-\bz^*\|^2\\
&\overset{\eqref{agradient-*}}{=}&  {\lambda} \|\bz-\bz^*\|^2.\end{array}
\end{eqnarray}
Thus $\bz^*$ is a unique local minimizer of (\ref{SLR-L2}).

For the case $\|\bz^*\|_0<s$, the condition \eqref{agradient-*} implies $\nabla f(\bz^*)=0$ due to $[\bz^*]^{\downarrow}_s=0$. Then \eqref{stongly-con-1} is true for any $\|\bz\|_0\leq s$. So $\bz^*$ is a unique global minimizer of the problem (\ref{SLR-L2}).

iii) ~ Let $\bz^*$ be a  $\tau$-stationary point for some $\tau >0$. Then  \eqref{alphaequation} and the definition of the projection $\Pi_s$ imply that
$$\|\bz^*-(\bz^*-\tau \nabla f(\bz^*))\|^2 \leq \|\bz-(\bz^*-\tau \nabla f(\bz^*))\|^2.$$
This leads to
\begin{equation*} 
2\langle \nabla f(\bz^*), \ \bz-\bz^* \rangle \geq -(1/\tau)\|\bz-\bz^*\|^2.
\end{equation*}
The above condition together with the inequality in \eqref{stongly-con-1} derives 
 \begin{eqnarray*} 
2f(\bz) \geq   2f(\bz^*)+ ({\lambda-1/\tau})\|\bz-\bz^*\|^2,
\end{eqnarray*}
which shows the unique global optimality of $\bz^*$ if $\tau>1/\lambda$. \qed

{We note that the necessary optimality condition in Theorem \ref{relation1} i) is directly adopted  from  \cite[Theorem 2.2]{Beck13} or \cite[Theorem 5.3]{beck15} where $\mathbb{B}=\mathbb{R}^p$. However, we also establish  the sufficient optimality conditions, see Theorem \ref{relation1} ii) and iii).} The above  relationships show that  a  $\tau$-stationary point is at least a unique locally optimal solution to the problem \eqref{SLR-L2}. This allows us to focus on a  $\tau$-stationary point itself to pursue a `good' solution. Therefore, we define a set
\begin{eqnarray}\label{Tu}
 \Sigma_s(\bz):=\{~\alpha\in [p]:~ |\alpha|=s,~  |z_i| \geq |z_j|, ~\forall~ i\in \alpha,~ j\in\overline{\alpha}~ \}. 
 \end{eqnarray}
Each element $\alpha$ in $\Sigma_s(\bz)$ coincides the  indices of the first $s$ largest (in absolute) components of $\bz$. Note that $\Sigma_s(\bz)$ may have multiple elements. For instance, $\bz=(3, -2, 2, 1, 0)^T$, $\Sigma_3(\bz)=\{\{1,2,3\}\}$ and $\Sigma_2(\bz)=\{\{1,2\},\{1,3\}\}.$ The notation allows us to rewrite  $\Pi_s(\bz)$ as follows
\begin{eqnarray*}
\Pi_s(\bz)= \Big\{\left(\bz_\alpha^T~0\right)^T:~ \alpha\in\Sigma_s(\bz)\Big\}.
\end{eqnarray*}
Then a point satisfying (\ref{alphaequation}) can be interpreted as that there is an $\alpha\in\Sigma_s(\bz-\tau \nabla f(\bz))$ satisfying $\bz_\alpha= (\bz-\tau \nabla f(\bz))_\alpha$ and $
\bz_{\overline{\alpha}}=0$,  which is equivalent to
\begin{eqnarray}\label{equation}
(\nabla f(\bz))_\alpha=0,~~~~
\bz_{\overline{\alpha}}=0.
\end{eqnarray}
Therefore, to find a  $\tau$-stationary point of (\ref{SLR-L2}), one can seek for a solution to the equation system \eqref{equation}.  This is summarized into the following theorem.
 \begin{theorem}\label{tau-F} A point $\bz$ is a $\tau$-stationary point of \eqref{SLR-L2} if and only if there is an $\alpha\in\Sigma_s(\bz-\tau \nabla f(\bz))$ satisfying \eqref{equation}.
\end{theorem}

\section{ {Convergence Analysis}}\label{Section3}
In this section, we turn our attention to solve the equations \eqref{equation} to {pursue} a $\tau$-stationary point of the problem \eqref{SLR-L2}, at least a unique local minimizer.

Given a   point $\bz^k$, for notational convenience, let
\begin{eqnarray*}
 H^k:=\nabla^2 f(\bz^k),~~\bg^k:=\nabla  f(\bz^k).\end{eqnarray*}
\subsection{The framework}  
Suppose we have a   point $\bz^k$ computed already. Then we can pick an index set $\alpha$ from $\Sigma_s(\bz^k-\tau {\bg}^k)$. For such a fixed index set $\alpha$, we apply Newton step on the equations \eqref{equation} just once to derive the Newton direction by
\begin{eqnarray*}
\begin{pmatrix}
H^k_{\alpha\alpha} &H^k_{\alpha\overline\alpha} \\
         0&I_{p-s}
\end{pmatrix} \bd^k= -\left(\begin{array}{c}
        \bg^k_{\alpha} \\
         \bz^k_{\overline\alpha}
       \end{array}\right)=: -\theta_\alpha^k,\end{eqnarray*}
where $I_{p-s}$ is the $(p-s)$ order identity matrix. One can calculate $\bd^k$ by
\begin{eqnarray}\label{newton-dir}
 \begin{cases}
  H^k_{\alpha\alpha} \bd^k_{\alpha} &=  H^k_{\alpha\overline\alpha}\bz^k_{\overline\alpha} -\bg^k_{\alpha},\\
~~~~~\bd^k_{\overline\alpha}&= -
         \bz^k_{\overline\alpha}.
       \end{cases} \end{eqnarray} 
Since $f$ is strongly convex, $H^k$ is non-singular and so are its any principal sub-matrices, i.e., $H^k_{\alpha\alpha}$ is invertible for any $k$ and any $\alpha$. Now we have the direction. If the full Newton step size is adopted, i.e., $\bz^{k+1}=\bz^{k}+\bd^{k}$, then \eqref{newton-dir} implies 
$$\|\bz^{k+1}\|_0=\|\bz^{k+1}_{\alpha} \|_0\leq |{\alpha} |=s.$$
Therefore, the updated point $\bz^{k+1}$ is feasible to the problem \eqref{SLR-L2}.  However, the full Newton step size generally does not guarantee the descent property of the objective {function}, that is, $f(\bz^{k+1})\leq f(\bz^{k})$ can not be ensured. To overcome such a drawback, we exploit the following operator
\begin{eqnarray}\label{z-sigma}
\bz^k(\sigma):=\begin{pmatrix}
        \bz^k_{\alpha} + \sigma \bd^k_{\alpha} \\
         \bz^k_{\overline\alpha} + \bd^k_{\overline\alpha}
       \end{pmatrix} =\begin{pmatrix}
        \bz^k_{\alpha} + \sigma \bd^k_{\alpha} \\
         0
       \end{pmatrix}. \end{eqnarray} 
For some carefully chosen $\sigma_k$, we set $\bz^{k+1}:=\bz^k(\sigma_k)$. Then  we will show that in this way,  $\bz^{k+1}$ is not only always feasible to the problem \eqref{SLR-L2} but also satisfies the descent property (see Lemma \ref{descent-direction}). Now we summarize the whole framework of Newton method in \Cref{Alg-NSLR}.

\begin{algorithm} 
	\caption{{\tt NSLR}: Newton method for SLR}
	\begin{algorithmic}[1] \label{Alg-NSLR}
		\STATE \textbf{Initialize}  $\bz^0,\tau>0, c\in(0,1)$. Set $k :=0$.
	\WHILE{the halting condition is violated}	
	    \STATE Pick  an $\alpha\in\Sigma_s(\bz^k-\tau \bg^k)$. \label{Step2}
	 \STATE  Update   $ \bd^k$ by \eqref{newton-dir}. \label{Step3}
	  \STATE   Find the smallest non-negative integer $r$ such that 
	  
	\vspace{-3mm}   
\parbox{0.9\textwidth}{\begin{eqnarray*}
2f(\bz^{k}(c^r))\leq 2f(\bz^k)+ c^{r} \langle \bg^k,  \bd^{k} \rangle.
\end{eqnarray*}}
 \STATE    Set  $\sigma_k:=c^{r}, \bz^{k+1}:=\bz^k(\sigma_k)$ and $k:=k+1$.
\ENDWHILE
\RETURN $\bz^k$
	\end{algorithmic}
\end{algorithm}
\begin{remark}\label{rem1}
Regarding   \NSLR, we have some comments.
\begin{itemize}
  \item [i)] To pick an $\alpha\in\Sigma_s(\bz^k-\tau \bg^k)$, only $s$ largest elements (in absolute) are selected, which enables us to use  a MATLAB built-in function \texttt{mink}. The computational complexity is about $\mathcal{O}(p+s\log s)$.

  \item [ ii)] Updating $\bd^{k}$ by (\ref{newton-dir}) involves two main calculations $H^k_{\alpha\alpha}$ with $|\alpha|=s$ and its inverse. Their computational complexities are $\mathcal{O}(sn+s^2n)$ and $\mathcal{O}(s^3)$. So, the whole complexity   is $\mathcal{O}(s^3+s^2n)$, which means the computation is quite fast if $\max\{s,n\}\ll p$.
    \item [{{iii)}}]  For a halting condition, we will calculate the quantity $\|\theta_\alpha^k\|$. If $\|\theta_\alpha^k\|=0$, then $\bz^k$ satisfies \eqref{equation}. This means $\bz^k$ is a $\tau$-stationary point by Theorem \ref{tau-F}. Therefore, it makes sense to terminate  \NSLR\ if the quantity $\|\theta_\alpha^k\|$ is sufficiently small.     
\end{itemize}
\end{remark}
\subsection{Global and quadratic convergence}To derive the convergence properties, we denote some notation hereafter. Let $\beta$ be the index set related the previous iteration $\bz^{k-1}$ selected by
$$\beta\in\Sigma_s(\bz^{k-1}-\tau { \bg}^{k-1}).$$
Based on $\bz^{k}=\bz^{k-1}(\sigma_{k-1}), \bz^{k+1}=\bz^{k}(\sigma_{k})$ in \Cref{Alg-NSLR} and  the definition of $\bz^{k}(\sigma )$ in \eqref{z-sigma}, we must have
\begin{eqnarray}
\label{supp-beta}
\supp(\bz^{k}) \subseteq \beta,~~~\supp(\bz^{k+1}) \subseteq \alpha.
\end{eqnarray}  
Let ${\gamma}:=\beta\setminus {\alpha}$.  Then one can observe that
\begin{eqnarray}
\label{z-k-alpha-com}
 -\bd^{k}_{\overline{\alpha}}=\bz^{k}_{\overline{\alpha}}  = \left(
 \begin{array}{c}
 \bz^{k}_{\beta\cap\overline{\alpha}}\\
 0
 \end{array}
 \right)= \left(
 \begin{array}{c}
 \bz^{k}_{\beta\setminus {\alpha}}\\
 0
 \end{array}
 \right)= \left(
 \begin{array}{c}
 \bz^{k}_{\gamma}\\
 0
 \end{array}
 \right).
\end{eqnarray}  
This gives rise to the following properties
\begin{eqnarray}\label{facts1}
\|\bd^k_{\overline\alpha}\|=\|\bd^k_{\gamma}\|=\|\bz^k_{\gamma}\|=\|\bz^k_{\overline\alpha}\|,~~ \langle \bg^k_{\overline\alpha}, \bd^k_{\overline\alpha}\rangle =\langle \bg^k_{\gamma}, \bd^k_{\gamma}\rangle. 
 \end{eqnarray}
 Based on these, we have the following results.
\begin{lemma}\label{properties}
Let $\{\bz^k\}$ be the sequence generated by \Cref{Alg-NSLR}. We have the following properties:
\begin{eqnarray}  \label{bounded-H}    \arraycolsep=1.4pt\def\arraystretch{1.25}
 \begin{array}{lcl}
2 \langle\bd^{k}_{\alpha },  \bg^{k}_{\alpha} \rangle &\leq&  - 2\lambda\| \bd^{k}_{ \alpha }\|^2+L \|\bd^{k}_{ {\overline\alpha} }\|^2 , \\
2 \langle\bd^{k}_{ {\overline\alpha}},  \bg^k_{ {\overline\alpha}} \rangle   &\leq& \tau L^2\| \bd^k_{\alpha} \|^2 +(\tau L^2-{1}/{\tau})\|  \bd^k_{\overline\alpha}\|^2.\end{array}
\end{eqnarray} 
\end{lemma}
{\bf Proof}
It follows from \eqref{newton-dir} that  
\begin{eqnarray}\label{fact-0}\bg^k_{\alpha}&=& -H^k_{\alpha\alpha}\bd^k_{\alpha}+ H^k_{\alpha\overline\alpha}\bz^k_{\overline\alpha} 
\overset{\eqref{newton-dir}}{=}
-H^k_{\alpha: }\bd^k.
 \end{eqnarray}
Direct calculation yields the following chain of equations,
       \begin{eqnarray*}
        \arraycolsep=1.4pt\def\arraystretch{1.25}
 \begin{array}{lcl}
       &&\langle \bd^{k},  H^k\bd^{k}\rangle-\langle   \bd^{k}_{ {\overline\alpha} }, H_{{\overline\alpha}{\overline\alpha} }^k\bd^{k}_{  {\overline\alpha} }\rangle\\
       &=& \langle\bd^{k}_{\alpha },  H_{\alpha\alpha }^k\bd^{k}_{\alpha}  \rangle + 2  \langle H^k_{\alpha{\overline\alpha}} \bd^{k}_{ {\overline\alpha} }, \bd^{k}_{\alpha }\rangle \\
        &=& 2 \langle H_{\alpha\alpha }^k\bd^{k}_{\alpha} +H^k_{\alpha{\overline\alpha}} \bd^{k}_{ {\overline\alpha} }, \bd^{k}_{\alpha }\rangle-\langle\bd^{k}_{\alpha },  H_{\alpha\alpha }^k\bd^{k}_{\alpha}  \rangle  \\
&\overset{\eqref{z-k-alpha-com}}{=}& 2 \langle H_{\alpha\alpha }^k\bd^{k}_{\alpha} -H^k_{\alpha{\overline\alpha}} \bz^{k}_{ {\overline\alpha} }, \bd^{k}_{\alpha }\rangle-\langle\bd^{k}_{\alpha },  H_{\alpha\alpha }^k\bd^{k}_{\alpha}  \rangle     \\
 &\overset{\eqref{fact-0}}{=}&- 2 \langle\bd^{k}_{\alpha },  \bg^k_{\alpha} \rangle - \langle H^k_{\alpha\alpha} \bd^{k}_{ \alpha }, \bd^{k}_{\alpha }\rangle,
 \end{array}
\end{eqnarray*}
which leads to the truth
 \begin{eqnarray*}       
  \arraycolsep=1.4pt\def\arraystretch{1.25}
 \begin{array}{lll}
2 \langle\bd^{k}_{\alpha },  \bg^{k}_{\alpha} \rangle&=& \langle   \bd^{k}_{ {\overline\alpha} }, H_{{\overline\alpha}{\overline\alpha} }^k\bd^{k}_{  {\overline\alpha} }\rangle - \langle H^k_{\alpha\alpha} \bd^{k}_{ \alpha }, \bd^{k}_{\alpha }\rangle -\langle \bd^{k},  H^k\bd^{k}\rangle\nonumber\\
&\leq&  L \|\bd^{k}_{ {\overline\alpha} }\|^2 - \lambda\| \bd^{k}_{ \alpha }\|^2- \lambda\| \bd^{k}_{\alpha\cup {\overline\alpha} }\|^2\nonumber\\
&\leq&  L \|\bd^{k}_{ {\overline\alpha} }\|^2 - 2\lambda\| \bd^{k}_{ \alpha }\|^2,
 \end{array}
\end{eqnarray*}
where the inequality is from the fact that $f$ being strongly convex with the constant $\lambda$ and strongly smooth with the constant $L$ so as to satisfy 
\begin{eqnarray}\label{eig}\lambda\leq\lambda_{\min}(H^k)\leq \lambda_{\max}(H^k)\leq L.\end{eqnarray}
For the part $\overline\alpha$, since $\alpha\in\Sigma_s(\bz^k-\tau \bg^k)$, the definition of $\alpha\in\Sigma_s$ in \eqref{Tu} implies
\begin{eqnarray*} 
 \forall~ i\in \alpha,~~ | z_i^k-\tau g^k_i| \geq |z_j^k-\tau g^k_j|, ~~ \forall~ j\in\overline{\alpha}.  
 \end{eqnarray*}
Now for $i\in \alpha\setminus \beta =: \eta$, we have $z_i^k=0$ due to \eqref{supp-beta}. Then the above condition and ${\gamma}\subseteq \overline{\alpha}$ result in
\begin{eqnarray*}
 \forall~ i\in \eta,~~ |\tau g^k_i| \geq |z_j^k-\tau g^k_j|, ~~ \forall~ j\in {\gamma},  
 \end{eqnarray*}
  which together with $\bz_{{\overline\alpha}}^k=- \bd_{{\overline\alpha}}^k$ from \eqref{newton-dir} and $|\eta|=|\alpha|-|\alpha \cap \beta|=s-|\alpha \cap \beta|=|{\beta}|-|\alpha \cap \beta|=|
 {\gamma}|$ leads to $$\|\tau \bg^k_{\eta}\|^2 \geq \|\bz_{{\gamma}}^k-\tau \bg^k_{{\gamma}}\|^2\overset{\eqref{z-k-alpha-com}}{=} \|\bd_{\gamma}^k+\tau \bg^k_{\gamma}\|^2.$$
The above condition allows us to derive that
 \begin{eqnarray*}
          \arraycolsep=1.4pt\def\arraystretch{1.25}
 \begin{array}{lcl}
 2 \langle\bd^{k}_{ {\overline\alpha}},  \bg^k_{ {\overline\alpha}} \rangle  &\overset{\eqref{facts1}}{=}& 2 \langle\bd^{k}_{{\gamma}},  \bg^{k}_{{\gamma}} \rangle\leq \tau \| \bg^k_{\eta}\|^2-\tau \| \bg^k_{\gamma}\|^2 -({1}/{\tau})\|  \bd^k_{\gamma}\|^2\nonumber\\
& {\leq}& \tau  \| \bg^k_{\alpha}\|^2 -({1}/{\tau})\|  \bd^k_{\gamma}\|^2 \nonumber\\
&\leq& \tau L^2\| \bd^k \|^2 -({1}/{\tau})\|  \bd^k_{\overline\alpha}\|^2\\
&=& \tau L^2\| \bd^k_{\alpha} \|^2 +(\tau L^2-{1}/{\tau})\|  \bd^k_{\overline\alpha}\|^2,
 \end{array}
\end{eqnarray*}
where the last inequality is owing to \eqref{facts1} and
  \begin{eqnarray}\label{nabla-T-f-0}
 \arraycolsep=1.4pt\def\arraystretch{1.25}
 \begin{array}{lcl}
 \|\bg^k_{\alpha}\|^2  &\overset{\eqref{fact-0}}{=}&  \|H^k_{\alpha: }  \bd^k\|^2 
~\leq ~ \|H^k_{\alpha: }  \bd^k\|^2 +\|H^k_{\overline\alpha: }  \bd^k\|^2 \\
& =&\|H^k  \bd^k\|^2   ~\overset{\eqref{eig}}{\leq}~ L^2\| \bd^k \|^2,
 \end{array} \end{eqnarray}
 showing the desired results.  \qed

The first result  below shows  that the Newton direction $\bd^k$ is a descent direction and the Amijio-type step size $\sigma_k$ always exists and is away from zero.
\begin{lemma}[Descent property] \label{descent-direction}
Let $\{\bz^k\}$ be the sequence generated by \Cref{Alg-NSLR} and  
\begin{eqnarray}\label{tau}
0<\tau<\overline{\tau}:=\min \left\{~ \frac{2\lambda}{L^2},~  \frac{c\lambda^2}{4L^3}~\right\}. 
 \end{eqnarray}
Denote $C:=\min\left\{ {1}/{\tau}-L-\tau L^2,~ 2\lambda-\tau L^2\right\}>0$. Then we have
 \begin{eqnarray}\label{decreasing-direction-gd}
 2\langle \bg^k, \bd^{k}\rangle  \leq -  C\| \bd^{k}\|^2.
 \end{eqnarray}
Moreover,  for any $c \overline{\sigma} \leq\sigma \leq  \overline{\sigma}:=\lambda/(2L)$, it holds 
\begin{eqnarray}\label{Armijo-sigma}
2f(\bz^{k}(\sigma))\leq 2f(\bz^k)+  \sigma \langle \bg^k,  \bd^{k} \rangle.
\end{eqnarray}
This indicates $\inf_{k\geq 0} \sigma_k \geq  c\overline{\sigma}>0$.
\end{lemma}
{\bf Proof}~ Denote two parameters
$$C_1:=\frac{1}{\tau}-L-\tau L^2,~~C_2:= 2\lambda-\tau L^2,~~C=\min\{C_1,C_2\}.$$
  It is easy to see that $C_2>0$ by \eqref{tau} and $C_1>0$ due to $$0<\tau<\overline{\tau}\leq  {c\lambda^2}/({4L^3})= ({c\lambda^2}/{L^2}) /({4L} )\leq {1}/({4L})$$ and $0<c, \lambda/L < 1$. Overall, $C>0$. Then by \eqref{bounded-H}, we have
\begin{eqnarray*}
\arraycolsep=1.4pt\def\arraystretch{1.25}
 \begin{array}{lll}
 2 \langle\bd^{k},  \bg^{k} \rangle &=& 2 \langle\bd^{k}_{\alpha },  \bg^{k}_{\alpha} \rangle+
2 \langle\bd^{k}_{\overline\alpha },  \bg^k_{\overline\alpha} \rangle\\
&{\leq}&    
-  C_1\| \bd^{k}_{ \alpha }\|^2 -   C_2  \|  \bd^k_{{\overline\alpha}}\|^2 \\
&\leq& -   C\| \bd^{k}_{ \alpha }\|^2-   C   \|  \bd^k_{{\overline\alpha}}\|^2\\
& = & -   C\| \bd^{k} \|^2. 
\end{array}\end{eqnarray*} 
Note from \eqref{z-sigma} that
\begin{eqnarray}\label{xk-alpha}
\bz^k(\sigma)-\bz^k=\left(\begin{array}{c}
         \sigma \bd^k_{\alpha} \\
         \bd^k_{\overline\alpha}
       \end{array}\right). \end{eqnarray} 
The strong  smoothness of $f$ with the constant $L$  yields
\begin{eqnarray*} 
\arraycolsep=1.4pt\def\arraystretch{1.25}
 \begin{array}{lll}
 &&2f(\bz^{k}(\sigma))-2f(\bz^k)-\sigma \langle \bg^k, \bd^k \rangle\\
&\leq &   2\langle \bg^k,   \bz^{k}(\sigma)-\bz^k   \rangle +L\|\bz^{k}(\sigma)-\bz^k\|^2-\sigma \langle \bg^k, \bd^k \rangle\nonumber\\
&= &  \sigma \langle \bg^k_\alpha,   \bd^k_\alpha \rangle+(2-\sigma)\langle \bg^k_{\overline\alpha},     \bd^k_{\overline\alpha} \rangle+L\sigma^2\|\bd^{k}_\alpha \|^2+L \|\bd^{k}_{\overline\alpha} \|^2 \nonumber\\
&=:&F.\nonumber
\end{array}\end{eqnarray*}
We next show $F\leq0$. It follows from \eqref{bounded-H} that
\begin{eqnarray*}
\arraycolsep=1.4pt\def\arraystretch{1.25}
 \begin{array}{lll}
 2F &{\leq}&  \Big(- 2\sigma\lambda + (2-\sigma) \tau L^2 \Big)\|\bd^{k}_\alpha \|^2\\
&+&\Big(\sigma L + (2-\sigma) (\tau L^2 - {1}/{\tau})\Big) \|\bd^{k}_{ {\overline\alpha} }\|^2+2\sigma^2L\|\bd^{k}_\alpha \|^2+2L \|\bd^{k}_{{\overline\alpha}} \|^2\\
 &=& \Big(- 2\sigma\lambda + (2-\sigma) \tau L^2 + 2\sigma^2L\Big)\|\bd^{k}_\alpha \|^2 \\
&+& \Big(\sigma L + (2-\sigma)  (\tau L^2 - {1}/{\tau} ) +2L\Big) \|\bd^{k}_{ {\overline\alpha} }\|^2.
\end{array}\end{eqnarray*}
One can check that, by $0<\tau<\overline{\tau}\leq {c\lambda^2}/({4L^3})$, it follows
\begin{eqnarray*}
\arraycolsep=1.4pt\def\arraystretch{1.25}
 \begin{array}{lll}
  &&- 2\sigma\lambda + (2-\sigma) \tau L^2 + 2\sigma^2L\\ 
&\leq&- 2\sigma\lambda + (2-\sigma) {c\lambda^2}/({4L}) + 2\sigma^2L\\
&=&2L\sigma^2 - \Big( 2\lambda +   {c\lambda^2}/({4L})\Big)  \sigma +  {c\lambda^2}/({2L})
\leq0,
\end{array}\end{eqnarray*}
where the last inequality is true if $ c\lambda/(2L) \leq\sigma \leq \lambda/(2L)$. For the $\overline{\alpha}$ part,
\begin{eqnarray*}
\arraycolsep=1.4pt\def\arraystretch{1.25}
 \begin{array}{lll}
   && \sigma L + (2-\sigma) \left(\tau L^2 - {1}/{\tau}\right) +2L \\ 
&\leq&\sigma L + (2-\sigma) L +({\sigma-2})/{\tau}  +2L \\
&\leq& 4L -  {1}/{\tau}  \leq 0,
\end{array}\end{eqnarray*}
where the above three inequalities used the facts that
\begin{itemize}
\item[(a)] $\tau L^2\leq {c\lambda^2}/({4L}) \leq \lambda (c/4)(\lambda/{L}) \leq L$;
\item[(b)] $\sigma\leq 1$ due to $\sigma_k=c^r$ and $c\in(0,1)$ in \Cref{Alg-NSLR}; 
\item[(c)] $\tau\leq {c\lambda^2}/({4L^3})\leq 1/({4L})$.
\end{itemize}
  Therefore, $F\leq 0$, displaying \eqref{Armijo-sigma}. Then the Armijo step-size rule indicates that $\inf_{k\geq 0} \sigma_k \geq  {c\lambda}/({2L})>0$. \qed
  
Now we are ready to display the main results of the method including the global convergence to a $\tau$-stationary point, the support set identification, and the quadratic convergence rate.
\begin{theorem}[Global and quadratic convergence] \label{global-quadratic}
Let $\{\bz^k\}$ be the sequence generated by \Cref{Alg-NSLR} and $0<\tau<\overline{\tau}$. The following results hold.
\begin{itemize}
\item[i)] The whole sequence converges to a $\tau$-stationary point denoted by $\bz^*$, which is at least a local minimizer of  the problem \eqref{SLR-L2}.
\item[ii)]For sufficiently large $k$, the support sets of the sequence are identified by
\begin{eqnarray}
\label{supp-alpha}
\supp(\bz^*)=
\begin{cases}
\supp(\bz^k)=\alpha,&\| \bz^*\|_0=s,\\
\supp(\bz^k)\cap\alpha,&\| \bz^*\|_0<s.
\end{cases}
\end{eqnarray}
\item[iii)]The sequence converges to  $\bz^*$ quadratically, namely,
\begin{eqnarray*}
\|\bz^{k+1}-\bz^{*} \|\leq  M/(2\lambda)\|\bz^{k}-\bz^{*} \|^2.
\end{eqnarray*} \end{itemize}
\end{theorem}
{\bf Proof} i) Lemma \ref{descent-direction} shows that  $\sigma_k \geq  c \overline{\sigma} $ and results in
\begin{eqnarray*}
\arraycolsep=1.4pt\def\arraystretch{1.25}
 \begin{array}{lcl}
  2f(\bz^{k+1})
=2f(\bz^{k}(\sigma_k)) 
&\overset{\eqref{Armijo-sigma}}{\leq}& 2f(\bz^k)+ \sigma_k \langle \bg^k,  \bd^{k} \rangle\\
&\overset{\eqref{decreasing-direction-gd}}{\leq}& 2f(\bz^k)- \sigma_k C \| \bd^k\|^2\\
&\leq& 2f(\bz^k)-  \overline{\sigma}  cC \| \bd^k\|^2.
\end{array}\end{eqnarray*}
Then it follows from the above inequality that
\begin{eqnarray*}
\arraycolsep=1.4pt\def\arraystretch{1.25}
 \begin{array}{lll}
   \sum^{\infty}_{k=0} \overline{\sigma}  cC \| \bd^k\|^2  
& \leq&
 \sum^{\infty}_{k=0}\Big(2 f(\bz^k)-2f(\bz^{k+1})\Big)\\
 &= & 2f(\bz^0)-\lim _{k\rightarrow +\infty}2f(\bz^k)\\
 &\leq& 2f(\bz^0),
\end{array}\end{eqnarray*}
where the last inequality is due to $f$ is positive. 
Hence  $\| \bd^k\|\rightarrow0$, which
suffices to $\|\bz^{k+1}-\bz^k\|\rightarrow0$ since
$$ \|\bz^{k+1}-\bz^k\|^2\overset{(\ref{xk-alpha})}{=}  \sigma^2 \| \bd^k_{\alpha}\|^2+\| \bd^k_{\overline \alpha}\|^2\rightarrow0.$$
This also indicates $\| \bz^k_{\overline \alpha}\|^2= \| \bd^k_{\overline \alpha}\|^2\rightarrow0$ and by (\ref{nabla-T-f-0}) suffices to 
\begin{eqnarray}\label{g-alpha-0} \|\bg^k_{ \alpha} \|  \leq  L \|\bd^k\|\rightarrow0.\end{eqnarray} 
Let  $\{\bz^{k_\ell}\}$ be the convergent subsequence of $\{\bz^{k}\}$  that converges to $\bz^*$ and
$$\alpha_\ell\in\Sigma_s(\bz^{k_\ell}-\tau\bg^{k_\ell}),~~ \ell\geq 1.$$
Since there are only finitely many choices for $\alpha_\ell\subseteq [p]$, 
(re-subsequencing if necessary)  we may without loss of any generality
assume that the sequence of the index sets $\{ \alpha_\ell\}$ shares a same index set, denoted
as $\alpha_{\infty}$. That is 
\begin{eqnarray*}
   \alpha_{\ell}=\alpha_{\ell+1}=\ldots =\alpha_\infty.
\end{eqnarray*}
Now by letting $\bg^*:=\nabla f(\bz^*)$, one can show that 
\begin{eqnarray} \label{limit-z-tau-g-supp}
\|\bg^*_{\alpha_\infty}\|=\lim_{k_\ell\rightarrow\infty}\|\bg^{k_\ell}_{ \alpha_\infty}\|=\lim_{k_\ell\rightarrow\infty}\|\bg^{k_\ell}_{ \alpha_\ell}\| \overset{\eqref{g-alpha-0}}{ = }0.
 \end{eqnarray}
In addition, the definition of $\Sigma_s$ in \eqref{Tu} implies
\begin{eqnarray*} 
 \forall~ i\in \alpha_{\ell}=\alpha_\infty,~ | z_i^{k_\ell}-\tau g^{k_\ell}_i| \geq |z_j^{k_\ell}-\tau g^{k_\ell}_j|,~ \forall~ j\in\overline{\alpha}_{\ell} =\overline\alpha_\infty. 
 \end{eqnarray*}
Taking the limit of both sides of the above inequality yields 
\begin{eqnarray} \label{limit-z-tau-g}
 \forall~ i\in  \alpha_\infty,~~ | z_i^*| \geq |\tau g^*_j|, ~~ \forall~ j\in \overline\alpha_\infty. 
 \end{eqnarray}
Here, we used the facts that \eqref{limit-z-tau-g-supp} and $\|\bz^{k_\ell}_{ \overline\alpha_\infty}\|=\|\bz^{k_\ell}_{ \overline\alpha_\ell}\| \rightarrow 0$. Since  $ \bz^{k_\ell}  \rightarrow \bz^*$, we have $\supp(\bz^*)\subseteq \alpha_\infty$. 
\begin{itemize}
\item If $\|\bz^*\|_0=s$, then $\supp(\bz^*)= \alpha_\infty$, which by \eqref{limit-z-tau-g} derives $|\tau g^*_j| \leq [\bz^*]^{\downarrow}_s, j \notin \supp(\bz^*)$. This together with \eqref{limit-z-tau-g-supp} results in condition \eqref{agradient}.
\item If $\|\bz^*\|_0<s$,  then $\supp(\bz^*) \subset \alpha_\infty$, which by \eqref{limit-z-tau-g} delivers $|\tau g^*_j| \leq [\bz^*]^{\downarrow}_s=0, j \notin  \overline\alpha_\infty$. This by \eqref{limit-z-tau-g-supp} yields $\bg^*=0$, which also satisfies \eqref{agradient}.
  \end{itemize} 
Overall, both cases show $\bz^*$ is a $\tau$ stationary point of \eqref{SLR-L2}.  From Theorem \ref{relation1}, a $\tau$ stationary point of \eqref{SLR-L2} is a unique local minimizer, namely, $\bz^*$ is isolated. By \cite[Lemma 4.10]{more1983computing}, the whole sequence converges to $\bz^*$ because $\bz^*$ is isolated  and $\|\bz^{k+1}-\bz^k\|\rightarrow0$.  

ii) We proved that the whole sequence converges to $\bz^*$. Denote $\alpha_*:=\supp(\bz^*)$. If $\bz^*=0$, then the conclusion holds clearly due to $\alpha_*=\emptyset$. Consider $\bz^*\neq0$. For sufficiently large $k$, we must have $$\|\bz^k-\bz^*\|<\min_{i\in \alpha_*} |z^*_i|=: \delta.$$
If $\alpha_* \nsubseteq \supp(\bz^k)$, then there is an $i_0\in\alpha_*\setminus\supp(\bz^k)$ satisfying
$$\delta>\|\bz^k-\bz^*\|\geq | z_{i_0}^k-z_{i_0}^* |=|  z_{i_0}^* |\geq \delta,$$
which is a contradiction. Therefore, $\alpha_*\subseteq\supp(\bz^k)$. By (\ref{supp-beta}), we have $\supp(\bz^k)\subseteq \beta$, where $|\beta|=s$ by (\ref{Tu}).  Therefore, if $\|\bz^{*}\|_0=s$ then $\beta \equiv \supp(\bz^k)\equiv \alpha_*$ for any sufficiently large $k$. Particularly,
 $\alpha_*\equiv \supp(\bz^{k+1})\equiv \alpha$. If $\|\bz^{*}\|_0<s$ then $\alpha_*
 \subseteq \supp(\bz^k)$ and $\alpha_*
 \subseteq  \supp(\bz^{k+1})\subseteq \alpha$ from (\ref{supp-beta}). So \eqref{supp-alpha} is true.
 
 iii) For sufficiently large $k$, it follows from $\alpha_*\subseteq \alpha$ by ii) that $\bz^*_{\overline\alpha}=0$. Since $\bz^*$ is a $\tau$-stationary point,  \eqref{agradient} indicates $\bg^*=0$ if $\|\bz^*\|_0<s$ and  $\bg^*_{\alpha_*}=0$ if $\|\bz^*\|_0=s$. While for the latter case, there is $\alpha_*=\alpha$  by ii). Overall, we  have
\begin{eqnarray}\label{F0-s-1}
\bz^*_{\overline\alpha}=0, ~~  \bg^*_\alpha =0.\end{eqnarray}
For any $0\leq t \leq1$, denote $\bz(t):=\bz^* + t (\bz^k-\bz^*)$ and $H^k(t):=\nabla^2f(\bz^k(t))$. Then \eqref{hess-lip} derives
\begin{eqnarray}
\label{facts-0-4} \|H^k-H^k(t) \|  
 \leq   M\|\bz^k-\bz(t)\| 
 = (1-t)M\|\bz^k- \bz^*\|.
\end{eqnarray}
 Moreover, by Taylor expansion, one has
\begin{eqnarray}\label{facts-0-3}
 \bg ^k - \bg ^*= { \int}_0^1 H^k(t)(\bz ^k-\bz ^*)dt.
\end{eqnarray}
From (\ref{supp-beta}) and (\ref{F0-s-1}), we have $\bz^{k+1}_{ \overline \alpha }=\bz^*_{  \overline \alpha }=0$ and the following relations
 \begin{eqnarray}\label{facts-0-5-1}
\arraycolsep=1.4pt\def\arraystretch{1.25}
 \begin{array}{lll}
   \|\bz^{k+1}-\bz^*\|^2 
&=&\| \bz^{k+1}_{ \alpha }-\bz^*_{ \alpha }\|^2\overset{(\ref{xk-alpha})}{=}\| \bz^{k}_{ \alpha }-\bz^*_{ \alpha }+\sigma_k \bd^{k}_{ \alpha }\|^2\\
&=&\|(1-\sigma_k) (\bz^{k}_{ \alpha }-\bz^*_{ \alpha })+\sigma_k (\bz^{k}_{ \alpha }-\bz^*_{ \alpha }+ \bd^{k}_{ \alpha })\|^2\\
&\leq&(1-\sigma_k)\| \bz^{k}_{ \alpha }-\bz^*_{ \alpha }\|^2+\sigma_k \|\bz^{k}_{ \alpha }-\bz^*_{ \alpha }+\bd^{k}_{ \alpha }\|^2~~~~~~ \\
&{\leq}&(1-c \overline{\sigma} )\| \bz^{k} -\bz^* \|^2+  \overline{\sigma} \|\bz^{k}_{ \alpha }-\bz^*_{ \alpha }+\bd^{k}_{ \alpha }\|^2,
\end{array}\end{eqnarray}
where the first inequality is due to $\|\cdot\|^2$ being convex  and the last inequality is from Lemma \ref{descent-direction} that $c \overline{\sigma} \leq\sigma_k \leq  \overline{\sigma}$. For the second term in \eqref{facts-0-5-1}, we have
\allowdisplaybreaks \begin{eqnarray}\label{facts-7}
\arraycolsep=1.4pt\def\arraystretch{1.25}
 \begin{array}{lcl}
  \lambda \|\bz^{k}_{ \alpha }-\bz^*_{ \alpha }+\bd^{k}_{ \alpha }\|
&\overset{(\ref{newton-dir})}{=}&\lambda\| (H^k_{\alpha\alpha})^{-1} (H^k_{\alpha\overline\alpha}\bz^k_{\overline\alpha} -\bg^k_{\alpha}  )+\bz^{k}_{ \alpha }-\bz^*_{ \alpha }\|\\
&\leq&\lambda \| (H^k_{\alpha\alpha})^{-1}  \|\cdot\|  H^k_{\alpha\overline\alpha}\bz^k_{\overline\alpha} -\bg^k_{\alpha} +H^k_{\alpha\alpha}(\bz^{k}_{ \alpha }-\bz^*_{ \alpha })\|\\
&\overset{(\ref{eig})}{\leq}&\|H^k_{\alpha:} \bz^k -\bg^k_{\alpha }  -H^k_{\alpha\alpha}\bz^*_{ \alpha } \|\\
&\overset{(\ref{F0-s-1})}{=}& \|H^k_{\alpha:} \bz^k
- \bg^k_{\alpha}-H^k_{\alpha:} \bz^* +\bg^*_{\alpha} \|\\
&\overset{(\ref{facts-0-3})}{=}& \|H^k_{\alpha:}(\bz^{k}-\bz^*)-\int_0^1H^k_{\alpha:}(t)(\bz ^k-\bz )dt\|\\
&=&\|\int_0^1(H^k_{\alpha:}-H^k_{\alpha:}(t))(\bz ^k-\bz ^*)dt\|\\
&\leq& \int_0^1\|(H^k_{\alpha:}-H^k_{\alpha:}(t))(\bz ^k-\bz ^*)\|dt\\
&\leq& \int_0^1\|H^k-H^k(t)\|\cdot\|\bz ^k-\bz ^*\|dt\\
&\overset{(\ref{facts-0-4})}{\leq}&M\|\bz ^k-\bz ^*\|^2 \int_0^1(1-t)dt\\
 &=& (M/2)\|\bz ^k-\bz ^*\|^2,
\end{array}\end{eqnarray}
where the forth inequality used a fact that $$\|A_{\alpha:}\bz\|^2\leq\|A_{\alpha:}\bz\|^2+\|A_{\overline\alpha:}\bz\|^2=\|A\bz\|^2\leq \|A\|^2\|\bz\|^2.$$ It follows from $\bd^k_{\overline \alpha}=-\bz^k_{\overline \alpha}$ and (\ref{F0-s-1})  that $$\|\bz^k+\bd^k-\bz^*\| = \|\bz^k_{\alpha }+\bd^k_{\alpha }-\bz^*_{\alpha }\|,$$ leading to the following fact
\begin{eqnarray}\label{first-fact}  \frac{\|\bz^k+\bd^k-\bz^*\|}{\|\bz^k-\bz^*\|} = \frac{\|\bz^k_{\alpha }+\bd^k_{\alpha }-\bz^*_{\alpha }\|}{\|\bz^k-\bz^*\|} \overset{( \ref{facts-7})}{\leq} \frac{M\|\bz ^k-\bz ^*\|^2}{2\lambda\|\bz^k-\bz^*\|}\rightarrow 0. \end{eqnarray}
Now we have three facts: (\ref{first-fact}), $\bz^k\rightarrow\bz^*$ from i), and
$ \langle \bg^{k}, \bd^{k}\rangle  \leq -  C\| \bd^{k}\|^2 $ from Lemma \ref{descent-direction}. They and \cite[Theorem 3.3]{facchinei1995minimization}
 allow us to claim that  eventually the step size $\sigma_k$ determined by the Armijo rule is 1,
 namely, $\sigma_k=1$. Then it follows from (\ref{facts-0-5-1}) that
  \begin{eqnarray*}
\label{facts-8}  \|\bz^{k+1}-\bz^*\|^2 &\overset{( \ref{facts-0-5-1})}{\leq}&(1-\sigma_k)\| \bz^{k}_{ \alpha }-\bz^*_{ \alpha }\|^2+\sigma_k \|\bz^{k}_{ \alpha }-\bz^*_{ \alpha }+\bd^{k}_{ \alpha }\|^2\nonumber\\
&=&  \|\bz^{k}_{ \alpha }-\bz^*_{ \alpha }+\bd^{k}_{ \alpha }\|^2
\overset{( \ref{facts-7})}{\leq} M^2/(2\lambda)^2\|\bz ^k-\bz ^*\|^4,
\end{eqnarray*}
delivering the quadratic convergence property of the sequence. \qed
\section{Numerical Experiments}\label{Section5}
In this part, we will conduct extensive numerical experiments of  \NSLR\footnote{Available at \url{https://github.com/ShenglongZhou/NSLR}} by using MATLAB (R2017b) on a desktop of 8GB of memory and Inter Core i5 2.7Ghz CPU, against seven  leading solvers on  both synthetic and real datasets.
\subsection{Test examples}
We first adopt two types of randomly generated data: the one with  { the i.i.d. features $(\bx_{1},\bx_{2},\ldots,\bx_{n})\in\mathbb{R}^{p\times n}$} and the one with independent features with each of $\bx_{i}$ being generated by an autoregressive process \cite{hamilton1994time}. Then eight real datasets are taken into consideration to test the selected methods.
\begin{example}[Independent Data \cite{lu2013sparse,pan2017convergent}] \label{EX1} To generate data labels $\by\in\{0,1\}^n$, we first randomly divide $[n]$ into two parts  and  set $y_i=0$ for one part and $y_i=1$ for the other. Then the feature data is produced by
$$\bx_i=y_iv_i\textbf{1}+\bw_i,\ \ \ \ i\in[n]$$
with $\mathbb{R}\ni v_i\sim \mathcal{N} (0,1)$, $\mathbb{R}^p\ni\bw_i\sim \mathcal{N} (0,I_p)$, where $\mathcal{N} (0,{I}_p)$ is the normal distribution with mean zero and variance identity. Here, $\textbf{1}$ is the vector with all entries being ones.   Since the  sparse parameter $\bz^*\in\mathbb{R}^p$ is unknown, different $s(<n)$ will be tested to {pursue} a sparse solution.
\end{example}

\begin{example}[Correlated Data \cite{agarwal2010fast, bahmani2013greedy}]\label{EX2}
The sparse parameter $\bz^*\in\mathbb{R}^p$  has $s$ nonzero entries drawn independently from the standard Gaussian distribution. Each data sample $\bx_i=(x_{i1},x_{i2},\ldots,x_{ip})^T, i\in[n]$ is an independent instance of the random vector
generated by an autoregressive process  \cite{hamilton1994time} determined by
$$x_{i(j+1)}=\rho x_{ij}+\sqrt{1-\rho^2}v_{ij},~~ j \in[p-1]$$
with $x_{i1}\sim \mathcal{N} (0,1)$, $v_{ij}\sim \mathcal{N} (0,1)$ and $ \rho\in[0,1]$ being the correlation parameter. The data labels $y_i\in\{0,1\}$
are then drawn randomly according to the Bernoulli distribution with  the  conditional probability \eqref{LR-cond-prob}. 
\end{example}

\begin{table}[!th]
\caption{Details of eight real datasets. \label{Detail-datasets} } \vspace{-10mm}
\renewcommand{\arraystretch}{0.85}\addtolength{\tabcolsep}{4pt}
\begin{center}
\begin{tabular}{lrrrr}\\\hline
 &  &  & Training & Testing\\  
Data name&$n$ &$p$ & $m_1$ & $m_2$\\ \hline
\texttt{arcene}&100 &10,000&100&0\\
\texttt{colon-cancer}&62 &2,000&62&0\\
\texttt{news20.binary}&19,996&1,355,191&19,996&0\\
\texttt{newsgroup}&11,314 &777,811 &11,314&0\\
\texttt{duke breast-cancer}&42	&7,129&38&4\\
\texttt{leukemia}&72		&7,129&38&34\\
\texttt{gisette}&7,000&5,000&6,000&1,000\\
\texttt{rcv1.binary}&	40,242&47,236&20,242&20,000\\\hline
\end{tabular}
\end{center}
\end{table}

\begin{example}[Real data]\label{EX3} Eight real data sets are taken into consideration. They are 
summarized in Table \ref{Detail-datasets}, where  \texttt{arcene} and \texttt{newsgroup} are taken from UCI repository\footnote{\url{http://archive.ics.uci.edu/ml/index.php}\label{1}} and {\tt glmnet} package\footnote{\url{https://web.stanford.edu/~hastie/glmnet\_matlab/}\label{3}}, and the rest of them are LIBSVM data\footnote{\url{https://www.csie.ntu.edu.tw/~cjlin/libsvmtools/datasets/}\label{2}}. Moreover, all datasets are feature-wisely scaled to $[-1,1]$. 
All $-1$s in the label classes $\by$ are replaced by 0. The sizes of training data and testing data are denoted by $m_1$ and $m_2$ respectively.

\end{example}
\subsection{Implementation}
In the model \eqref{SLR-L2}, we set $\lambda=10^{-5}/n$. As mentioned in Remark \ref{rem1}, we terminate the method if  $\|\theta^{k}\|<10^{-10}\sqrt p$ or $k>2000$. For the starting point and parameters of \NSLR,   we set $\bz^0=0$ and $c=0.5$. For the parameter $\tau$, Theorem \ref{global-quadratic} indicates   $0<\tau<\overline{\tau}$. While this is a sufficient condition. So it is unnecessary to choose a $\tau$ to satisfy the condition strictly, not to mention, $\overline{\tau}$ might be too small if we set a tiny $\lambda$. 

\begin{figure}[!th]
\centering
\begin{subfigure}{.33\textwidth}
  \centering
  \includegraphics[width=1.01\textwidth]{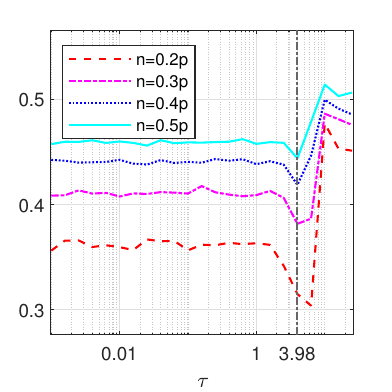}
\caption{$\ell(\bz)$}
\end{subfigure}
\begin{subfigure}{.33\textwidth}
  \centering
  \includegraphics[width=1.01\textwidth]{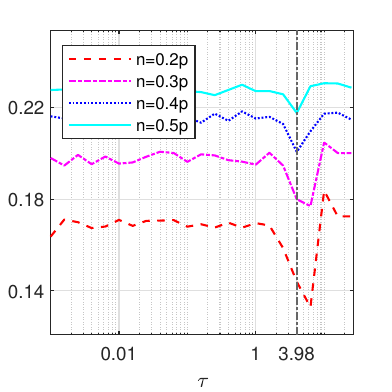}
   \caption{${\tt SER}$}
\end{subfigure}%
\begin{subfigure}{.33\textwidth}
  \centering
  \includegraphics[width=1.01\textwidth]{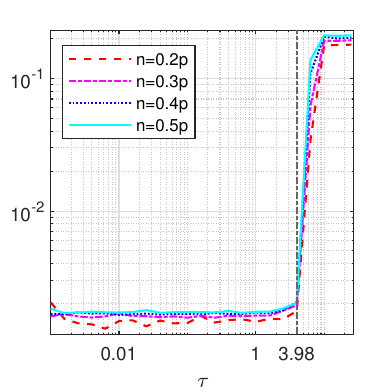}
   \caption{${\tt Time}$}
\end{subfigure} \\
\begin{subfigure}{.33\textwidth}
  \centering
  \includegraphics[width=1.01\textwidth]{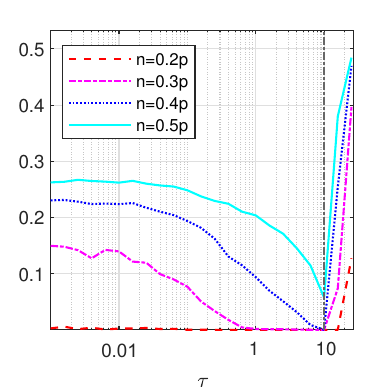}
\caption{$\ell(\bz)$}
\end{subfigure}
\begin{subfigure}{.33\textwidth}
  \centering
  \includegraphics[width=1.01\textwidth]{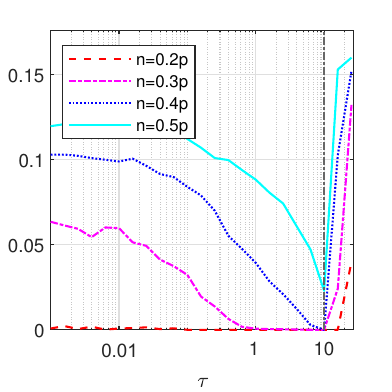}
   \caption{${\tt SER}$}
\end{subfigure}%
\begin{subfigure}{.33\textwidth}
  \centering
  \includegraphics[width=1.01\textwidth]{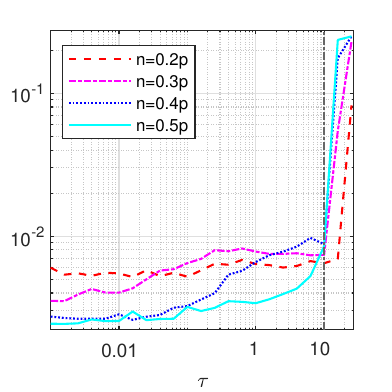}
   \caption{${\tt Time}$}
\end{subfigure}
\caption{Effect to $\tau$ of {\tt NSLR} for Example \ref{EX1} with  $p=500$, $s=5$ (the above three sub-figures) and  $s=25$ (the below three sub-figures).}
\label{fig:effect-tau}
\end{figure}

Alternative is to pick a proper $\tau$ by tuning it from a wide range of values. For instance, we tested \NSLR\ on a range of selections of $\tau=10^\varrho$ with $\varrho\in\{-3,-2.8,\ldots,1,1.2,1.4\}$ for solving Example \ref{EX1}. As reported in Fig. \ref{fig:effect-tau}, for $s=5$, \NSLR\ generated the best results when $\tau$ was around $3.98$, while for $s=25$, it delivered the best results when $\tau$ was around $10$. Therefore, for different scenarios, the best option $\tau$ may be varied, which indicates the manual selection of  $\tau$  is necessary to achieve a better performance. This apparently would incur expensive computational costs.

However,  the empirical numerical experience have demonstrated that adaptively updating parameters during the process is a practicable strategy.   Hence, we start  $\tau$ with a fixed scalar $\tau_0=15$ and  update $\tau_{k+1}=0.75\tau_k$ if $k$ is a multiple of 10 and $\|\theta^{k}\|> 1/k$, and $\tau_{k+1}=\tau_k $ otherwise. In the sequel, we adopt this strategy for \NSLR\ and the numerical comparisons with other leading solvers will show the superior performance of our method under such a strategy.

\subsection{Benchmark methods}
Since there is an impressive body of methods that have been developed  to address the sparse logistic regression,  we only focus on those programmed by Matlab. Solvers with codes being online unavailable or being written by other languages, such as R and C, are not selected for comparisons. We thus choose 7 solvers mentioned in  \Cref{methods-SLR}, which should be enough to make comprehensive comparisons. We summarize them into the following table. 
\begin{table}[H]
\caption{Benchmark Methods.}\vspace{-10mm} 
\renewcommand{\arraystretch}{1.0}\addtolength{\tabcolsep}{-2pt}
\begin{center}
\begin{tabular}{cccc}\\\hline
Models &(\ref{SLR-p}) with convex $\phi_\nu$   & (\ref{SLR-p}) with non-convex $\phi_\nu$ &(\ref{SLR-L2})\\\hline
first-order&\texttt{SLEP}&\texttt{APG},~\texttt{GIST}&\texttt{GraSP}\\
second-order&\texttt{IRLS-LARS}&\texttt{--}&\texttt{NTGP}, \texttt{GPGN}\\\hline
\end{tabular}\end{center}
\end{table}
\vspace{-5mm}

For \texttt{SLEP}, we use it to solve (\ref{SLR-p}) with $\phi_{\nu}(\bz)=\nu_1\|\bz\|^2_2+\nu_2\|\bz\|_1$, whilst \texttt{IRLS-LARS} aims to solve the case of  $\nu=0$.   \texttt{APG} and \texttt{GIST} are taken to solve the capped $\ell_1$  logistic regression with $\phi_{\nu}(\bz) = \nu_3\min(|\bz_i|,\nu_4)$. We only use non-monotonous version of \texttt{APG}  since its numerical performance was better than that of the  monotonous version  \cite{li2015accelerated}. Note that methods that aim at solving the model (\ref{SLR-p}) involve a penalty parameter $\nu$, whilst those tackling (\ref{SLR-L2}) need the sparsity level $s$. To make results comparable, we adjust their default parameters $\nu$ for each method to guarantee the generated solution $\bz$ satisfying $\|\bz\|_0\leq p/2$. We will report three indicators:
 $( \ell(\bz),~\ser,~\time )$
to illustrate the performance of methods, where \time\ (in seconds) is the CPU time, $\bz$ is the solution obtained by each method and \ser\ is the sign error rate defined by
$$\ser:=\frac{1}{n}\sum_{i=1}^m\Big|y_i-\mathrm{sign}\left(\langle \bx_i,\bz\rangle_+\right)\Big|.$$
Here $\mathrm{sign}(a_+)$ is the sign of the projection of $a$ onto a non-negative space, namely, it returns $1$ if $a>0$ and $0$ otherwise.

 \subsection{Numerical comparison}
 We now report the performance of eight methods on the above three examples. To avoid randomness, we report average results over 10-time independent trails for Examples \ref{EX1} and \ref{EX2} since they involve in randomly generated data.

\begin{figure}[!th]
\centering
\begin{subfigure}{.33\textwidth}
  \centering
  \includegraphics[width=1.01\textwidth]{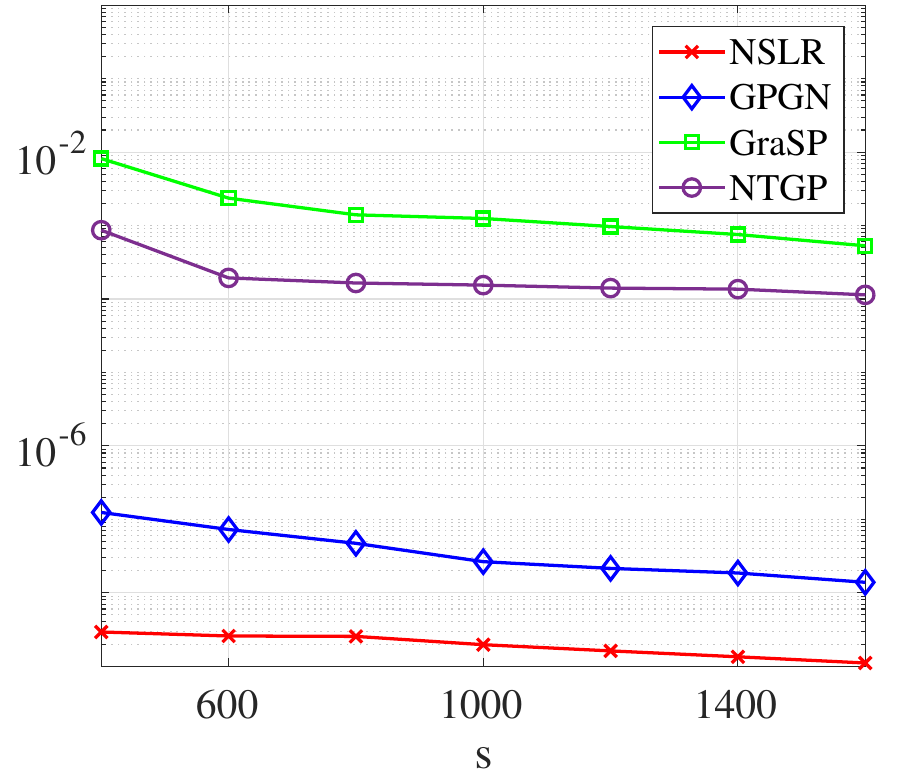}
  \caption{$\ell(\bz)$}
\end{subfigure}%
\centering
\begin{subfigure}{.33\textwidth}
  \centering
  \includegraphics[width=1.01\textwidth]{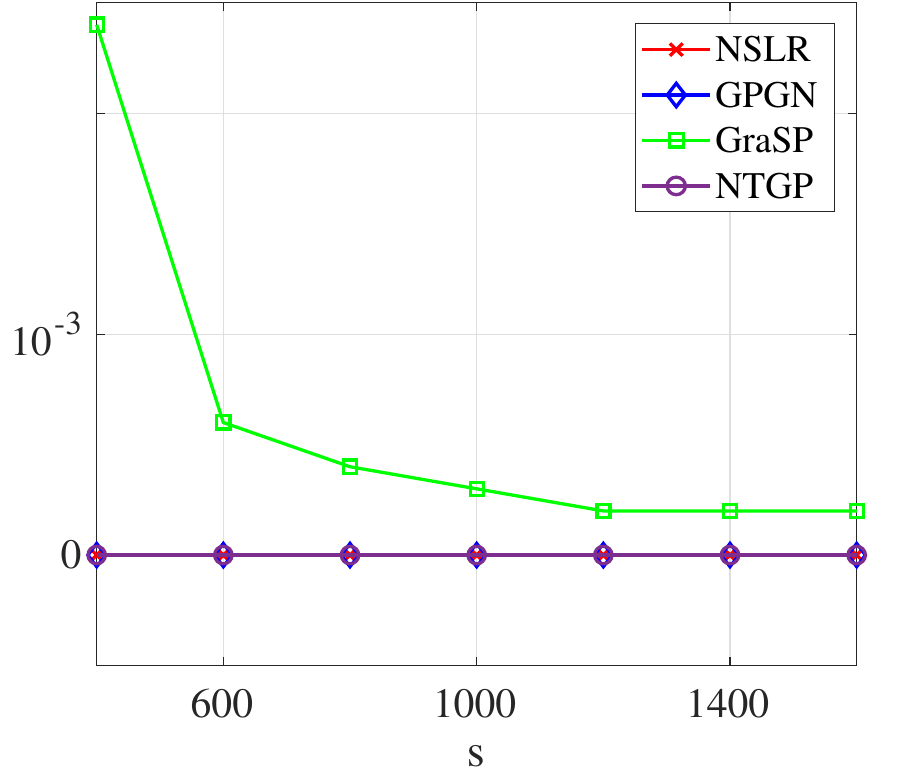}
   \caption{${\tt SER}$}
\end{subfigure}%
\begin{subfigure}{.33\textwidth}
  \centering
  \includegraphics[width=1.01\textwidth]{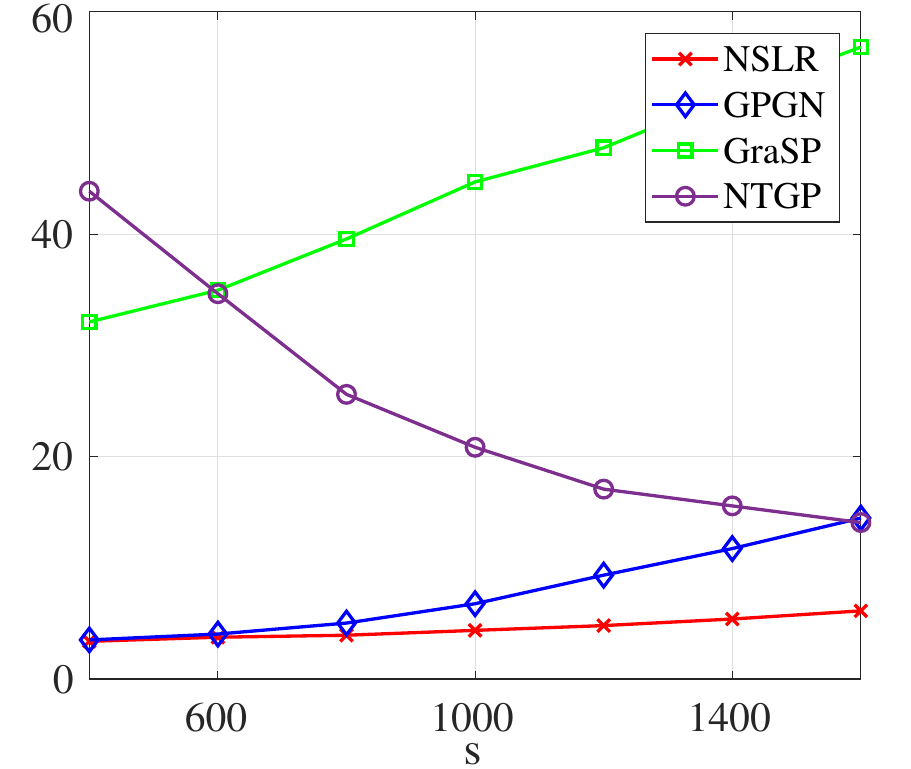}
   \caption{${\tt Time}$}
\end{subfigure}
\caption{Comparison of four methods for Example \ref{EX1} with $p = 10000, n=0.2p$}
\label{fig1}
\end{figure}

\textbf{(a) Comparison on Example \ref{EX1}.} To observe the influence of the sparsity level $s$ on four greedy methods: \NSLR, \GPGN, \GraSP\ and \NTGP,
we fix $p=10000, n=p/5$ and alter $s\in\{400,600,\ldots,1600\}.$ As demonstrated in Fig. \ref{fig1},  \NSLR\ outperforms others in terms of the lowest  $\ell(\bz)$ and \ser\ and the shortest time, followed by \GPGN. By contrast, \GraSP\ always performs the worst results, which means this first-order method is not competitive when against the other three methods, three second-order methods.

To observe the influence of the ratio of the sample size $n$ and the number of features $p$ on all eight methods, we fix $p=20000, \ s=0.1p$  and vary $n/p\in\{0.1,0.2,\ldots,0.7\}$. Apart from recording the four indicators, we also report the number of non-zeros of the solution $\bz$ generated by each method. Here,  for the \LARS, we stop it when $s$ variables are selected and their default stopping conditions are met since \LARS\ only adds one variable at each iteration (see \cite{huang2018constructive}). We set $\nu_1=10^{-1},\nu_2=10^{-2}$ for  \SLEP, $\nu_3=10^{-2},\ \nu_4=10^{-4}$ for \APG\ and $\nu_3=10^{-3}{\tt abs(randn)}, \ \nu_4=10^{-5}{\tt abs(randn)}$ for   \GIST. 

\begin{figure}[!th]
\centering
\begin{subfigure}{.33\textwidth}
  \centering
  \includegraphics[width=1.04\textwidth]{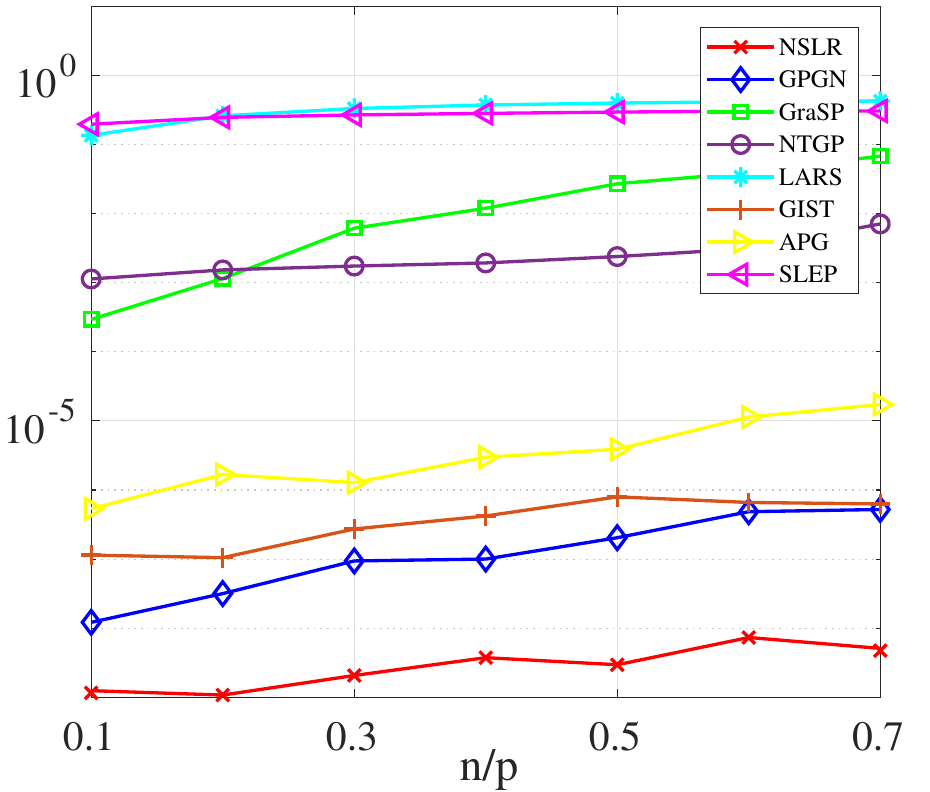}
\caption{$\ell(\bz)$}
\end{subfigure}
\centering
\begin{subfigure}{.33\textwidth}
  \centering
  \includegraphics[width=1.04\textwidth]{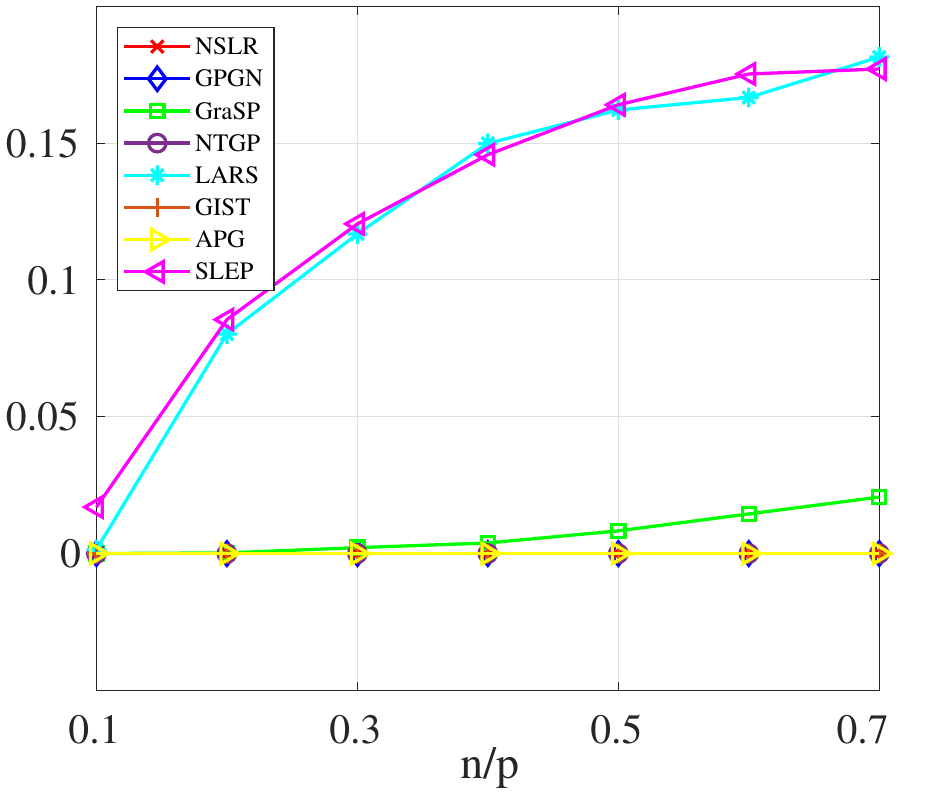}
   \caption{${\tt SER}$}
\end{subfigure}%
\begin{subfigure}{.33\textwidth}
  \centering
  \includegraphics[width=1.04\textwidth]{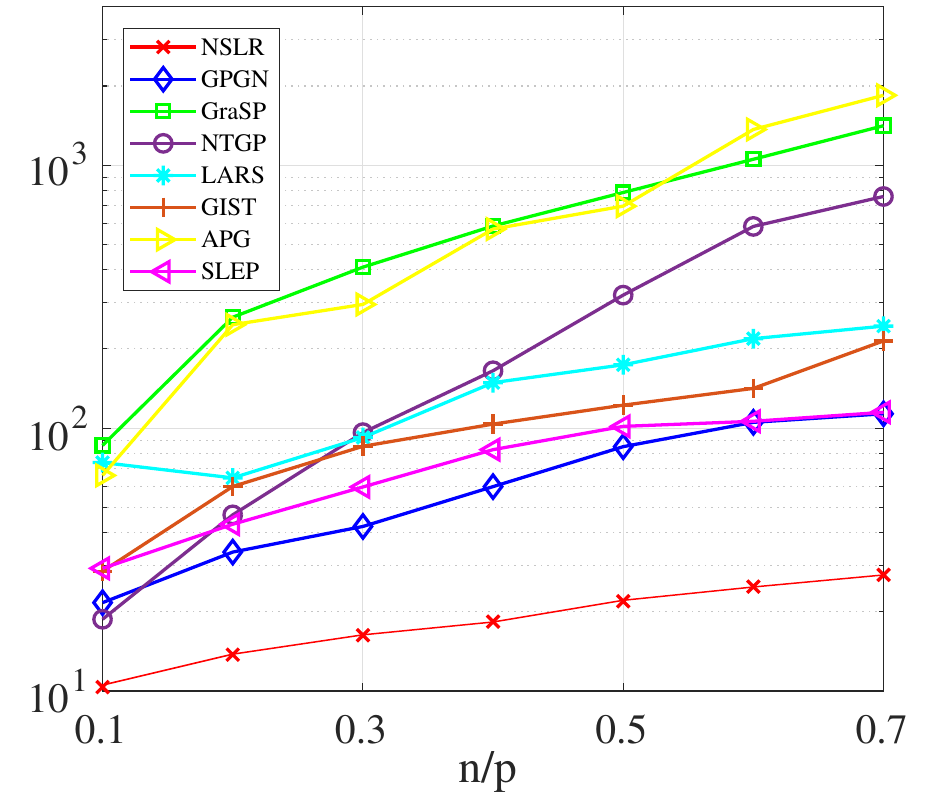}
   \caption{${\tt Time}$}
\end{subfigure}
\caption{Comparison of eight methods for Example \ref{EX1} with $p = 20000, s = 0.1p.$}
\label{fig2}
\end{figure}

As presented in Fig. \ref{fig2}, in terms of $\ell(\bz)$ and \ser, again \NSLR\ performs the best results, followed by  \GPGN\ and   \GIST. It is {obvious} that \LARS\ and \SLEP\ produce undesirable results compared with other methods. For the computational time, \NSLR\ runs the fastest, while \GraSP\  and \APG\ run relatively slow with over 1000 seconds when $n/p\geq 0.6$. Table \ref{tab0} shows the sparsity levels $\|\bz\|_0$ only in \LARS\ is lower than our \NSLR. This is because \LARS\ fails to recover the support and vanishes when $s = 500$ in this numerical experiment (this phenomenon had also been observed in \cite{huang2018constructive,garg2009gradient}.)

\begin{table} [!th]
\caption{Sparsity levels $\|\bz\|_0$ of eight methods for Example \ref{EX1} with $p = 20000, s = 0.1p$. \label{tab0}} 
\vspace{-3mm}
\renewcommand{\arraystretch}{0.8}\addtolength{\tabcolsep}{3pt}\begin{center}
\begin{tabular} {c rrrrrrr}\hline
$n/p$ &$0.1$ &$0.2$ &0.3 &0.4&0.5&0.6&0.7\\\hline
 \texttt{NSLR, GPGN} &\multirow{ 2}{*}{2000}	&\multirow{ 2}{*}{2000}	&\multirow{ 2}{*}{2000}&\multirow{ 2}{*}{2000}&\multirow{ 2}{*}{2000} &\multirow{ 2}{*}{2000} &\multirow{ 2}{*}{2000}\\
  \texttt{GraSP, NTGP} &&&& \\
 \verb"LARS" &500	&500	&500 &500 &500 &500 &500\\
\verb"GIST"	&4403 	&4309 	&5274 	&7832 	&7913 	&8614 	&8904 	\\
\verb"APG"	&6138 	&5857 	&5720 	&6170 	&5574 	&5048 	&5043 	\\
\verb"SLEP"	&2076 	&2534 	&2980 	&3235 	&3498 	&3596 	&3873 	\\\hline
\end{tabular}\end{center}
\end{table}

\begin{figure}[!th]
\centering
\begin{subfigure}{.33\textwidth}
  \centering
  \includegraphics[width=1.01\textwidth]{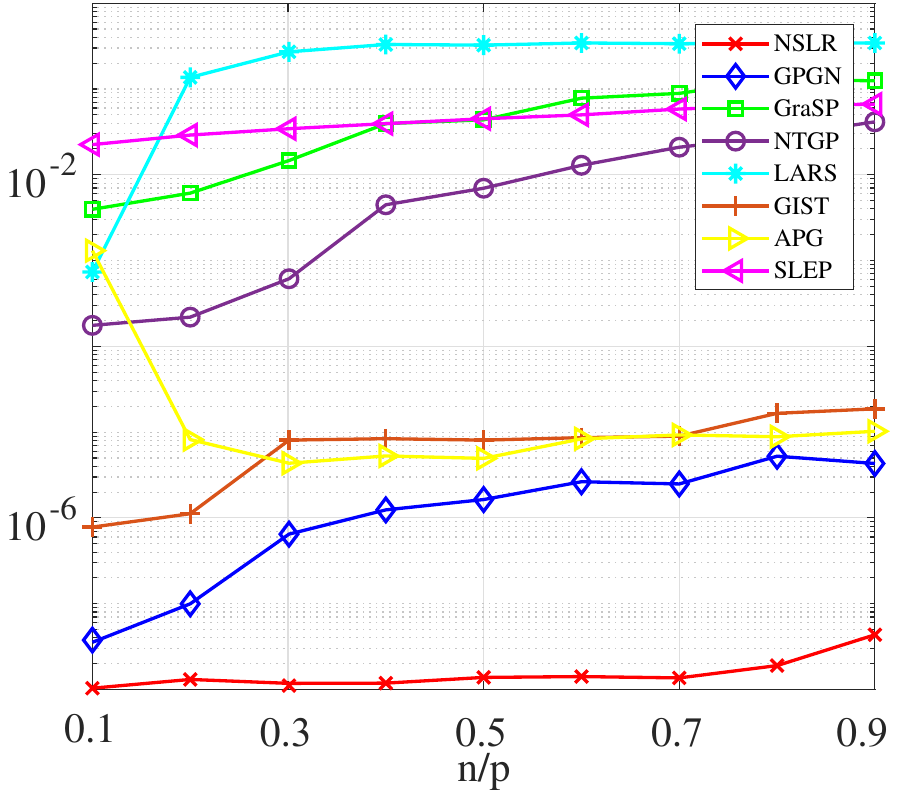}
  \small \texttt{ (a) $\rho=0$}
\end{subfigure}%
\begin{subfigure}{.33\textwidth}
  \centering
  \includegraphics[width=1.01\textwidth]{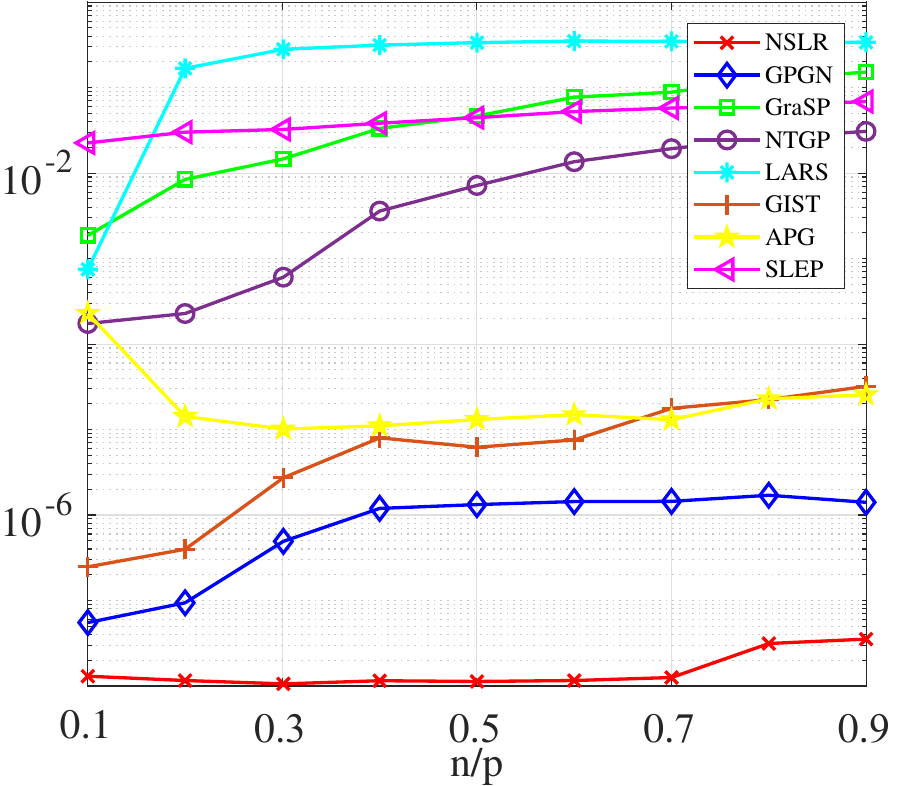}
 \small \texttt{(b) $\rho=1/3$}
\end{subfigure}
\centering
\begin{subfigure}{.33\textwidth}
  \centering
  \includegraphics[width=1.01\textwidth]{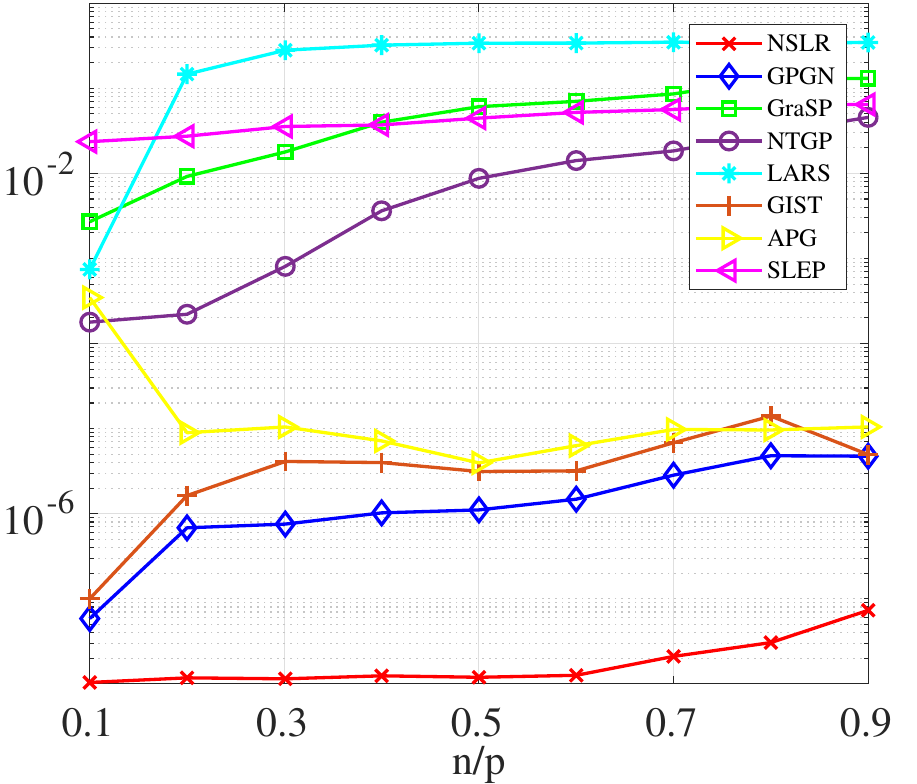}
   \small \texttt{ (c) $\rho=1/2$}
\end{subfigure}%
\caption{$\ell(\bz)$ obtained by eight methods for Example \ref{EX2} with $p = 1000, s = 0.1p.$}
\label{fig3}
\end{figure}

 \textbf{(b) Comparison on Example \ref{EX2}.}  To observe the influence of the correlation parameter $\rho$ on  eight methods, we set $p=1000, s=0.1p$ but choose $\rho\in\{0,1/3, 1/2\}$.  Fig. \ref{fig3} shows the average $\ell(\bz)$ gotten by eight methods for a wide range of the ratio $n/p\in\{0.1,0.2,\ldots,0.9\}$. Apparently, at three different values of $\rho$, \NSLR\ always performs stably best results. Moreover, the trends in these eight methods perform generally consistent, which indicates the correlation parameter has little influence on these methods. Therefore, we further fix $\rho= 1/2$ and observe the performance of eight methods under higher dimensions.

\begin{table}[!th] 
{  \caption{ Average results  for Example \ref{EX2}. \label{tab1}}\vspace{0mm}
{\renewcommand\baselinestretch{1}\selectfont
\renewcommand\tabcolsep{11pt}  
{\centering\begin{tabular}{l cccrccccr}\hline
  &\multicolumn{4}{c }{$s=0.05p, n=0.2p$}&&\multicolumn{4}{c}{$s=0.1p, n=0.2p$}  \\\cline{2-5}\cline{7-10}
&{$\ell(\bz)$} &{\ser} &Time&{$\|\bz\|_0$}&&{$\ell(\bz)$} &{\ser} &Time&{$\|\bz\|_0$} \\\cline{2-5}\cline{7-10}
&\multicolumn{8}{c }{$p=10000$} \\\hline
\verb"NSLR"	 &   3.2e-10 	 &   0.00e-0 	 &   0.436  	 &   500 &     	 &   1.1e-10 	  &   0.00e-0 	&    0.918  	&    1000  \\	
\verb"GPGN"	 &   1.09e-7	 &   0.00e-0	 &   4.192 	 &   500 &    	 &   6.84e-8	  &   0.00e-0	 &   6.962 	&   1000 	\\
\verb"GraSP"	 &   1.20e-2	 &   5.00e-3	 &   17.10 	 &   500 &    	 &   7.41e-3	  &   1.70e-3	 &   17.01 	 &   1000 	\\
\verb"NTGP"	 &   4.37e-3	 &   0.00e-0	 &   33.54 	 &   500 &    	 &   2.18e-3	  &   0.00e-0	 &   13.53 	 &   1000 	\\
\verb"LARS"	 &   1.43e-1	 &   1.25e-2	 &   27.69 	 &   500  &   	 &   2.85e-1	  &   4.90e-3	 &   71.80 	 &   1000	\\
\verb"GIST"	 &   4.81e-6	 &   0.00e-0	 &   12.54 	 &   2381 &    	 &   8.08e-6	  &   0.00e-0	 &   11.40 	 &   2974 	\\
\verb"APG"	 &   1.29e-4	 &   0.00e-0	 &   30.00 	 &   1407 &    	 &   1.25e-4	  &   0.00e-0	 &   27.82 	 &   5211 	\\
\verb"SLEP"	 &   1.75e-1	 &   2.00e-3	 &   6.560	 &   811 &    	 &   1.32e-1	  &   0.00e-0	 &   10.84 	 &   1019 	\\\hline
&\multicolumn{8}{c }{$p=20000$} \\\hline
\verb"NSLR"	 &    1.6e-10 	 &    0.00e-0 	 &    2.522  	 &   1000  &       	 &   5.4e-11 	  &   0.00e-0 	 &   4.939  	 &    2000  \\	
\verb"GPGN"	 &   8.48e-8	 &   0.00e-0	 &   16.47 	 &   1000 &    	 &   5.91e-8	  &   0.00e-0	 &   15.83 	 &   2000 	\\
\verb"GraSP"	 &   1.29e-2	 &   5.25e-3	 &   62.22 	 &   1000  &   	 &   5.00e-3	  &   1.50e-3	 &   92.69 	 &   2000 	\\
\verb"NTGP"	 &   5.43e-3	 &   0.00e-0	 &   134.6 	 &   1000 &    	 &   2.30e-3	  &   0.00e-0	 &   54.50 	 &   2000 	\\
\verb"LARS"	 &   3.96e-1	 &   5.43e-2	 &   107.4 	 &   1000 &   	 &   4.21e-1	  &   6.18e-2	 &   117.1 	 &   1000 	\\
\verb"GIST"	 &   7.20e-7	 &   0.00e-0	 &   33.52 	 &   1542 &    	 &   9.98e-7	  &   0.00e-0	 &   54.55 	 &   4137 	\\
\verb"APG"	 &   5.06e-5	 &   0.00e-0	 &   24.94 	 &   1849 &    	 &   6.04e-5	  &   0.00e-0	 &   69.46 	 &   4323 	\\
\verb"SLEP"	 &   1.88e-1	 &   3.25e-3	 &   22.67 	 &   1511 &    	 &   1.45e-1	  &   0.00e-0	 &   34.59 	 &   2005 	\\\hline
&\multicolumn{8}{c }{$p=30000$} \\\hline
\verb"NSLR"	 &    1.1e-10 	 &    0.00e-0 	 &   7.364 	  &   1500 &   	 &   3.8e-11	 &   0.00e-0	 &   12.79 	 &   3000	\\
\verb"GPGN"	 &   8.58e-8	 &   0.00e-0	 &   35.06 	 &   1500 &   	 &   4.09e-8	  &   0.00e-0	 &   51.49 	 &   3000	\\
\verb"GraSP"	 &   2.28e-2	 &   7.50e-3	 &   181.9 	 &   1500 &   	 &   8.96e-3	 &   2.83e-3	 &   208.1 	 &   3000	\\
\verb"NTGP"	 &   5.36e-3	 &   0.00e-0	 &   307.1 	 &   1500 &   	 &   2.32e-3	  &   0.00e-0	 &   125.1 	 &   3000	\\
\verb"LARS"	 &   4.59e-1	 &   1.03e-1	 &   205.7 	 &   1000 &   	 &   4.99e-1	 &   1.21e-1	 &   206.7 	 &   1000 	\\
\verb"GIST"	 &   2.76e-6	 &   0.00e-0	 &   121.1 	 &   6124 &   	 &   7.97e-7	  &   0.00e-0	 &   150.9 	 &   6278	\\
\verb"APG"	 &   9.40e-4	 &   1.67e-4	 &   206.7 	 &   7460 &   	 &   3.11e-4	  &   0.00e-0	 &   247.1 	 &   6450	\\
\verb"SLEP"	 &   1.96e-1	 &   2.33e-3	 &   54.20 	 &   2392 &   	 &   1.49e-1	  &   1.67e-4	 &   80.19 	 &   3009	\\\hline

\end{tabular}\par} }}
\end{table}
Now we  alter $p\in\{10000,20000,30000\}$ with $n=0.2p, s=0.05p$ or $s=0.1p$.
Here, we set $\nu_1=10^{-2},\nu_2=10^{-1}$ for \SLEP\, $\nu_3=10^{-2},\ \nu_4=5\times 10^{-4}$ for \APG\ and $\nu_3=5\times10^{-3}{\tt abs(randn)}, \ \nu_4=5\times10^{-5}{\tt abs(randn)}$ for \GIST. As reported in Table \ref{tab1}, for cases of $p=20000,30000$, \LARS\ basically fails to render desirable solutions due to the highest $\ell(\bz)$. It can be clearly seen that \NSLR\ always provides the best accuracies  with consuming shortest time. 


 \textbf{(c) Comparison on Example \ref{EX3}.} To observe the performance for above all eight methods on real data sets, we select eight real data sets with different dimensions. The highest dimension is up to millions (see \texttt{news20.binary}).  Table \ref{tab2} reports  results for eight methods on four datasets without testing data ($m_2=0$): \texttt{arcene}, \texttt{colon-cancer}, \texttt{news20.binary} and \texttt{newsgroup}. For the last two datasets, \LARS\ makes our desktop run  out of memory, thus  its results are omitted here.  Clearly, \NSLR\ is more efficient than others for all test instances. For example, \NSLR\ only uses 1.365 seconds for data \texttt{newsgroup} with $p=777811$ features and achieves the smallest logistic loss with the sparsest solution.

\begin{table}[!h]
{  \caption{ Results  for Example \ref{EX3} with $m_2=0$. \label{tab2}}
{\renewcommand\baselinestretch{1}\selectfont
\renewcommand\tabcolsep{11pt}  
{\centering\begin{tabular}{ l  cccrccccr }\hline
&{$\ell(\bz)$}  &{\ser} &{Time}&{$|| \bz ||_0$}& 
&{$\ell(\bz)$}  &{\ser} &{Time}&{$|| \bz ||_0$}\\\cline{2-5}\cline{7-10}
Data&\multicolumn{4}{c }{{\tt Arcene}}&&\multicolumn{4}{c}{{\tt colon-cancer}}  \\\hline
\verb"NSLR"	 &    4.57e-7 	 &    0.00e-0 	 &   0.141 	 &    60 &    	 &    1.90e-8 	 &    0.00e-0 	  &    0.082  	  &    20 	\\
 \verb"GPGN"	 &   1.84e-5	 &   0.00e-0	 &    0.583  	&   60 &   	 &   3.35e-5	 &   0.00e-0	 &   0.171 	&   20	\\
\verb"GraSP"	 &   1.88e-3	 &   0.00e-0	 &   1.233 	&   60 &   	 &   3.03e-1	  &   1.61e-2	 &   0.264 	&   20	\\
\verb"NTGP"	 &   1.09e-1	 &   1.00e-2	 &   2.812 	 &   60 &   	 &   4.51e-3	  &   0.00e-0	 &   0.192 	&   20	\\
\verb"LARS" 	 &   7.98e-2	 &   0.00e-0	 &   0.974 	&   60 &   	 &   2.84e-2	  &   0.00e-0	 &   0.231 	&   30	\\
\verb"GIST"	 &   3.99e-7	 &   0.00e-0	 &   5.644 	&   134 &   	 &   6.22e-7	  &   0.00e-0	 &   0.133 	&   41	\\
\verb"APG"	 &   1.84e-5	 &   0.00e-0	 &   3.341 	&   430 &   	 &   1.79e-7	  &   0.00e-0	 &   0.234 	&   20	\\
\verb"SLEP"	 &   1.05e-1	 &   0.00e-0	 &   3.663 	&   66 &   	 &   1.68e-1	  &   4.84e-2	 &   0.142 	&   23	\\\hline

Data&\multicolumn{4}{c }{ {\tt news20.binary}}  &&\multicolumn{4}{c}{{\tt  newsgroup}}\\\hline
\verb"NSLR"	 &    1.11e-2 	 &    3.50e-3 	 &    3.236  	  &    2500 &    	 &    1.46e-2 	 &    8.00e-3 	  &    1.365  	  &    3000 	\\
 \verb"GPGN"	 &   2.94e-2	 &   8.00e-3	 &   25.44 	 &   2500 &   	 &   5.17e-2	  &   1.20e-2	 &   30.71 	 &   3000	\\
\verb"GraSP"	 &   2.46e-2	 &   9.25e-3	 &   205.7 	 &   2500 &   	 &   5.08e-1	&   4.75e-2	 &   33.06 	 &   3000	\\
\verb"NTGP"	 &   7.94e-2	 &   6.55e-3	 &   93.20 	 &   2500 &   	 &   3.11e-1	  &   5.80e-2	 &   13.07 	 &   3000	\\
\verb"LARS"&   $--$	&   $--$	&   $--$	&  $--$ &	&   $--$	&   $--$	&   $--$	&  $--$	\\
\verb"GIST"	 &   3.18e-2	 &   6.60e-3	 &   27.57 	 &   4091 &   	 &   6.84e-2	  &   1.52e-2	 &   10.55 	 &   3017	\\
\verb"APG"	 &   4.18e-2	 &   1.24e-2	 &   37.15 	 &   3869 &   	 &   5.12e-2	  &   1.77e-2	 &   10.62 	 &   3257	\\
\verb"SLEP"	 &   1.49e-1	 &   3.34e-2	 &   43.70 	 &   5299 &   	 &   2.60e-1	  &   4.60e-2	 &   23.89 	 &   5138	\\\hline

\end{tabular}\par} }}
\end{table}

When all methods solve the datasets with testing data ($m_2>0$): \texttt{duke breast-cancer},
\texttt{leukemia}, \texttt{gisette} and \texttt{rcv1.binary}, the table is a little different.     Since the testing data is taken into consideration, we add two indicators to illustrate the performance of each method: $\ell(\bz)$-test and  \ser-test. Results are reported in Table \ref{tab3}, where $ \ell(\bz) $ and $ \ell(\bz) $-test denote the objective function value on training data and testing data respectively, and similar to \ser-train and  \ser-test. For  cases of \texttt{duke breast-cancer} and \texttt{leukemia}, \LARS\ stops when the maximum number of iterations reaches the $\min\{m_1,p\}$. Hence its produced $|| \bz ||_0$s are less than other methods.  We can see that \NSLR\ could guarantee a good performance on the testing data as well as the training data. 

\section{Conclusion}\label{Section6}
  Despite the NP-hardness of the sparsity constrained logistic regression (\ref{SLR-L2}),
we benefited from its nice properties of the objective function and introduced the   $\tau$-stationary point as an optimality condition. This can be converted to an equation system that makes the Newton method effective.  Since an order $s$ (which is far smaller than $p$)  principal sub-matrix of the whole Hessian of the objective function is taken into account in each step, the proposed method \NSLR\ has a relatively low computational complexity. The success in acquiring the global convergence stems from the realization of the Armijo-type line search in the method. What is more, the generated sequence also converges to a $\tau$-stationary point quadratically, which well justifies the outstanding performance of \NSLR\ theoretically. It is worth mentioning that we reasonably extended the classical Newton method for solving unconstrained and continuous problems to the sparsity constrained logistic regression. The numerical performance against several state-of-the-art methods demonstrated that \NSLR\ is remarkably efficient and competitive, especially in large scale settings.

We feel that the proposed method might be capable of solving the general strong convex optimization problems with the sparsity constraint. This deserves exploring in future.

\begin{table}[!th]
{  \caption{Results for Example \ref{EX3} with $m_2>0$. \label{tab3}}
{\renewcommand{\arraystretch}{.75}\addtolength{\tabcolsep}{11pt}
{\centering\begin{tabular}{ lcccccr}\hline
&{$\ell(\bz)$}  &{$ \ell(\bz) $-test} &{\ser-train}  &{\ser-test} &{Time}&{$||\bz||_0$}\\\hline
{Data} &\multicolumn{6}{c}{{\tt duke breast-cancer}}  \\\hline
\verb"NSLR"	 &  4.10e-9 	 &  2.45e-6 	 &  0.00e-0 	 &  0.00e-0 	 &  0.113  	 &  100	\\
 \verb"GPGN"	 &  1.31e-5	 &  2.55e-3	 &  0.00e-0	 &  0.00e-0	 &  0.442 	 &  100	 \\
\verb"GraSP"	 &  3.57e-3	 &  1.82e-3	 &  0.00e-0	 &  0.00e-0	 &  0.502 	 &  100	 \\
\verb"NTGP"	 &  1.21e-5	 &  1.20e-4	 &  0.00e-0	 &  0.00e-0	 &  0.641 	 &  100	\\
\verb"LARS" 	 &  1.01e-4	 &  1.21e-4	 &  0.00e-0	 &  0.00e-0	 &  0.613 	 &  37 	\\
\verb"GIST"	 &  6.64e-9	 &  1.60e-1	 &  0.00e-0	 &  0.00e-0	 &  0.544 	 &  614	\\
\verb"APG"	 &  6.35e-7	 &  2.68e-7	 &  0.00e-0	 &  0.00e-0	 &  0.742 	 &  136	\\
\verb"SLEP"	 &  1.93e-3	 &  1.04e-2	 &  0.00e-0	 &  0.00e-0	 &  0.434 	 &  203	 \\\hline

{Data} &\multicolumn{6}{c}{{\tt leukemia}}  \\\hline
\verb"NSLR"	 &  3.09e-6 	 &  7.22e-2	 &  0.00e-0 	 &  0.00e-0 	 &  0.113  	 &  150	\\
 \verb"GPGN"	 &  1.29e-5	 &  2.71e-0	 &  0.00e-0	 &  1.47e-1	 &  0.351 	 &  150	 \\
\verb"GraSP"	 &  3.40e-3	 &  5.08e-1	 &  0.00e-0	 &  1.47e-1	 &  0.562 	 &  150	 \\
\verb"NTGP"	 &  4.22e-4	 &  2.11e-1	 &  0.00e-0	 &  5.88e-2	 &  0.754 	 &  150	\\
\verb"LARS" 	 &  6.07e-4	 &  1.11e-1 	 &  0.00e-0	 &  8.82e-2	 &  1.052 	 &  37 	\\
\verb"GIST"	 &  5.68e-3	 &  1.82e-1	 &  0.00e-0	 &  1.18e-1	 &  0.423 	 &  295	\\
\verb"APG"	 &  1.05e-4	 &  3.63e-1	 &  0.00e-0	 &  8.82e-2	 &  0.542 	 &  1066	\\
\verb"SLEP"	 &  1.71e-1	 &  2.85e-1	 &  0.00e-0	 &  2.94e-2	 &  0.581 	 &  269	 \\\hline

{Data} &\multicolumn{6}{c}{{\tt gisette}}  \\\hline
\verb"NSLR"	&  2.88e-6       	&  6.51e-1 	 &  0.00e-0  	 &  4.30e-2 	 &  1.221 	 &  500 	\\
 \verb"GPGN"	 &  2.25e-4	 &  2.63e-1	 &  0.00e-0	 &  4.60e-2	 &  1.412 	 &  500	 \\
\verb"GraSP"	 &  1.51e-4	 &  2.30e-1	 &  0.00e-0	 &  4.82e-2	 &  2.433 	 &  500	 \\
\verb"NTGP"	 &  1.17e-3	 &  9.58e-1	 &  0.00e-0	 &  4.18e-2	 &  2.554 	 &  500	\\
\verb"LARS" 	 &  1.63e-1	 &  1.39e-0	 &  1.00e-3	 &  4.18e-2	 &  8.262 	 &  500	 \\
\verb"GIST"	 &  2.40e-4	 &  1.02e-0	 &  0.00e-0	 &  4.30e-2	 &  1.783 	 &  1303	\\
\verb"APG"	 &  3.72e-4	 &  1.31e-0	 &  0.00e-0	 &  4.90e-2	 &  1.641 	 &  907	\\
\verb"SLEP"	 &  1.92e-2	 &  8.30e-1	 &  0.00e-0	 &  4.97e-2	 &  2.042 	 &  1569	 \\\hline

Data &\multicolumn{6}{c}{{\tt rcv1.binary}}  \\\hline
\verb"NSLR"	 &  5.82e-2	 &  2.01e-1 	 &  1.95e-2	 &  5.45e-2 	 &  3.712  	 &  1000  	\\
 \verb"GPGN"	 &  2.90e-2	 &  2.35e-1	 &  8.05e-3	 &  5.58e-2	 &  6.811 	 &  1000	 \\
\verb"GraSP"	 &  3.09e-1	 &  1.98e-0	 &  4.15e-2	 &  9.51e-2	 &  21.72 	 &  1000	 \\
\verb"NTGP"	 &  7.48e-2	 &  1.37e-1	 &  1.06e-2	 &  4.64e-2	 &  4.471 	 &  1000	\\
\verb"LARS" 	 &  2.27e-1	 &  2.61e-1	 &  5.13e-1	 &  5.46e-1	 &  33.53 	 &  1000	 \\
\verb"GIST"	 &  3.44e-2	 &  1.41e-1	 &  8.20e-3	 &  4.85e-2	 &  4.571 	 &  1545	\\
\verb"APG"	 &  2.51e-2 	 &  2.78e-1	 &  7.31e-3  	 &  6.00e-2	 &  7.571 	 &  1537	\\
\verb"SLEP"	 &  1.24e-1	 &  1.65e-1	 &  3.02e-2	 &  5.23e-2	 &  8.453	 &  3527	\\\hline
\end{tabular}\par} }}
\end{table}

\section*{Acknowledgements}
This work is supported by the National
Natural Science Foundation of China (11971052) and Beijing Natural Science Foundation (Z190002). We particularly thank the referee and the Principal Editor who offered us valuable suggestions to improve this paper greatly.

\bibliography{ref}

\begin{thebibliography}{10}
\expandafter\ifx\csname url\endcsname\relax
  \def\url#1{\texttt{#1}}\fi
\expandafter\ifx\csname urlprefix\endcsname\relax\def\urlprefix{URL }\fi
\expandafter\ifx\csname href\endcsname\relax
  \def\href#1#2{#2} \def\path#1{#1}\fi

\bibitem{tibshirani1996regression}
R.~Tibshirani, Regression shrinkage and selection via the {L}asso, J. R. Stat.
  Soc. Series B Stat. Methodol. (1996) 267--288.

\bibitem{bahmani2013greedy}
S.~Bahmani, B.~Raj, P.~T. Boufounos, Greedy sparsity-constrained optimization,
  J. Mach. Learn. Res. 14~(Mar) (2013) 807--841.

\bibitem{plan2013robust}
Y.~Plan, R.~Vershynin, Robust 1-bit compressed sensing and sparse logistic
  regression: A convex programming approach, IEEE Trans. Inf. Theory 59~(1)
  (2013) 482--494.

\bibitem{Wang}
R.~Wang, N.~Xiu, C.~Zhang, Greedy projected gradient-newton method for
  large-scale sparse logistic regression, IEEE Trans. Neural Netw. Learn. Syst.
  31~(2) (2019) 527--538.

\bibitem{beck15}
A.~Beck, N.~Hallak, On the minimization over sparse symmetric sets:
  projections, optimality conditions, and algorithms, Math. Oper. Res. 41~(1)
  (2015) 196--223.

\bibitem{pan2015solutions}
L.~Pan, N.~Xiu, S.~Zhou, On solutions of sparsity constrained optimization, J.
  Oper. Res. Soc. Chn. 3~(4) (2015) 421--439.

\bibitem{Beck13}
A.~Beck, Y.~Eldar, Sparsity constrained nonlinear optimization: Optimality
  conditions and algorithms, SIAM J. Optim. 23~(3) (2013) 1480--1509.

\bibitem{hastie2017extended}
T.~Hastie, R.~Tibshirani, R.~Tibshirani, Extended comparisons of best subset
  selection, forward stepwise selection, and the {L}asso, arXiv preprint
  arXiv:1707.08692 (2017).

\bibitem{mazumder2017subset}
R.~Mazumder, P.~Radchenko, A.~Dedieu, Subset selection with shrinkage: Sparse
  linear modeling when the {SNR} is low, arXiv preprint arXiv:1708.03288
  (2017).

\bibitem{hazimeh2018fast}
H.~Hazimeh, R.~Mazumder, Fast best subset selection: Coordinate descent and
  local combinatorial optimization algorithms, Oper. Res. 68~(5) (2020)
  1517--1537.

\bibitem{xie2018ccp}
W.~Xie, X.~Deng, The {CCP} selector: Scalable algorithms for sparse ridge
  regression from chance-constrained programming, arXiv preprint
  arXiv:1806.03756 (2018).

\bibitem{Pang2019}
T.~Pang, F.~Nie, J.~Han, X.~Li, Efficient feature selection via $\ell
  _{2,0}$-norm constrained sparse regression, IEEE Trans. Knowl. Data Eng.
  31~(5) (2019) 880--893.

\bibitem{figueiredo2003adaptive}
M.~Figueiredo, Adaptive sparseness for supervised learning, IEEE Trans. Pattern
  Anal. Mach. Intell. 25~(9) (2003) 1150--1159.

\bibitem{krishnapuram2005sparse}
B.~Krishnapuram, L.~Carin, M.~Figueiredo, A.~Hartemink, Sparse multinomial
  logistic regression: Fast algorithms and generalization bounds, IEEE Trans.
  Pattern Anal. Mach. Intell. 27~(6) (2005) 957--968.

\bibitem{andrew2007scalable}
G.~Andrew, J.~Gao, Scalable training of $\ell_1$-regularized log-linear models,
  in: Proceedings of the 24th International Conference on Machine Learning,
  2007, pp. 33--40.

\bibitem{koh2007interior}
K.~Koh, S.~Kim, S.~Boyd, An interior-point method for large-scale
  $\ell_1$-regularized logistic regression, J. Mach. Learn. Res. 8~(Jul) (2007)
  1519--1555.

\bibitem{yu2010quasi}
J.~Yu, S.~Vishwanathan, S.~G{\"u}nter, N.~Schraudolph, A quasi-newton approach
  to nonsmooth convex optimization problems in machine learning, J. Mach.
  Learn. Res. 11~(Mar) (2010) 1145--1200.

\bibitem{shi2010fast}
J.~Shi, W.~Yin, S.~Osher, P.~Sajda, A fast hybrid algorithm for large-scale
  $\ell_1$-regularized logistic regression, J. Mach. Learn. Res. 11~(Feb)
  (2010) 713--741.

\bibitem{yuan2010comparison}
G.~Yuan, K.~Chang, C.~Hsieh, C.~Lin, A comparison of optimization methods and
  software for large-scale $\ell_1$-regularized linear classification, J. Mach.
  Learn. Res. 11~(Nov) (2010) 3183--3234.

\bibitem{liu2009slep}
J.~Liu, S.~Ji, J.~Ye, {SLEP}: Sparse learning with efficient projections,
  Arizona State University 6~(491) (2009) 7.

\bibitem{friedman2010regularization}
J.~Friedman, T.~Hastie, R.~Tibshirani, Regularization paths for generalized
  linear models via coordinate descent, J. Stat. Softw. 33~(1) (2010) 1--22.

\bibitem{yuan2012improved}
G.~Yuan, C.~Ho, C.~Lin, An improved {GLMNET} for $\ell_1$-regularized logistic
  regression, J. Mach. Learn. Res. 13~(Jun) (2012) 1999--2030.

\bibitem{liu2009large}
J.~Liu, J.~Chen, J.~Ye, Large-scale sparse logistic regression, in: Proceedings
  of the 15th ACM SIGKDD International Conference on Knowledge Discovery and
  Data Mining, 2009, pp. 547--556.

\bibitem{lee2006efficient}
S.~Lee, H.~Lee, P.~Abbeel, A.~Ng, Efficient $\ell_1$ regularized logistic
  regression, in: Proceedings of the 21st National Conference on Artificial
  Intelligence and the 18th Innovative Applications of Artificial Intelligence
  Conference, 2006, pp. 401--408.

\bibitem{efron2004least}
B.~Efron, T.~Hastie, I.~Johnstone, R.~Tibshirani, Least angle regression, Ann.
  Stat. 32~(2) (2004) 407--451.

\bibitem{fan2001variable}
J.~Fan, R.~Li, Variable selection via nonconcave penalized likelihood and its
  oracle properties, J. Am. Stat. Assoc. 96~(456) (2001) 1348--1360.

\bibitem{zou2008one}
H.~Zou, R.~Li, One-step sparse estimates in nonconcave penalized likelihood
  models, Ann. Stat. 36~(4) (2008) 1509--1533.

\bibitem{huang2009group}
J.~Huang, S.~Ma, H.~Xie, C.~Zhang, A group bridge approach for variable
  selection, Biometrika 96~(2) (2009) 339--355.

\bibitem{gong2013general}
P.~Gong, C.~Zhang, Z.~Lu, J.~Huang, J.~Ye, A general iterative shrinkage and
  thresholding algorithm for non-convex regularized optimization problems, in:
  30th International Conference on Machine Learning, 2013, pp. 37--45.

\bibitem{li2015accelerated}
H.~Li, Z.~Lin, Accelerated proximal gradient methods for nonconvex programming,
  in: Advances in Neural Information Processing Systems 28-Proceedings of the
  2015 Conference, 2015, pp. 379--387.

\bibitem{gong2015honor}
P.~Gong, J.~Ye, {HONOR}: Hybrid optimization for non-convex regularized
  problems, in: Advances in Neural Information Processing Systems, 2015, pp.
  415--423.

\bibitem{rakotomamonjy2016dc}
A.~Rakotomamonjy, R.~Flamary, G.~Gasso, Dc proximal newton for nonconvex
  optimization problems, IEEE Trans. Neural Netw. Learn. Syst. 27~(3) (2016)
  636--647.

\bibitem{needell2009cosamp}
D.~Needell, J.~Tropp, {CoSaMP}: Iterative signal recovery from incomplete and
  inaccurate samples, Appl. Comput. Harmon. Anal. 26~(3) (2009) 301--321.

\bibitem{lozano2011group}
A.~Lozano, G.~Swirszcz, N.~Abe, Group orthogonal matching pursuit for logistic
  regression, in: Proceedings of the 14th International Conference on
  Artificial Intelligence and Statistics, Vol.~15, 2011, pp. 452--460.

\bibitem{mallat1993matching}
S.~Mallat, Z.~Zhang, Matching pursuits with time-frequency dictionaries, IEEE
  Trans. Signal Process. 51~(12) (1993) 3397--3415.

\bibitem{lu2013sparse}
Z.~Lu, Y.~Zhang, Sparse approximation via penalty decomposition methods, SIAM
  J. Optim. 23~(4) (2013) 2448--2478.

\bibitem{yuan2014gradient}
X.~Yuan, P.~Li, T.~Zhang, Gradient hard thresholding pursuit for
  sparsity-constrained optimization, in: 31st International Conference on
  Machine Learning, 2014, pp. 1322--1330.

\bibitem{pan2017convergent}
L.~Pan, S.~Zhou, N.~Xiu, H.~Qi, A convergent iterative hard thresholding for
  nonnegative sparsity optimization, Pac. J. Optim. 13~(2) (2017) 325--353.

\bibitem{yuan2017newton}
X.~Yuan, Q.~Liu, Newton-type greedy selection methods for $\ell_0$-constrained
  minimization, IEEE Trans. Pattern Anal. Mach. Intell. 39~(12) (2017)
  2437--2450.

\bibitem{chen2017fast}
J.~Chen, Q.~Gu, Fast newton hard thresholding pursuit for sparsity constrained
  nonconvex optimization, in: Proceedings of the 23rd ACM SIGKDD International
  Conference on Knowledge Discovery and Data Mining, ACM, 2017, pp. 757--766.

\bibitem{more1983computing}
J.~Mor{\'e}, D.~Sorensen, Computing a trust region step, SIAM J. Sci. Comput.
  4~(3) (1983) 553--572.

\bibitem{facchinei1995minimization}
F.~Facchinei, Minimization of {SC1} functions and the maratos effect, Oper.
  Res. Lett. 17~(3) (1995) 131--138.

\bibitem{hamilton1994time}
J.~Hamilton, Time series analysis, Vol.~2, Princeton university press
  Princeton, NJ, 1994.

\bibitem{agarwal2010fast}
A.~Agarwal, S.~Negahban, M.~Wainwright, Fast global convergence rates of
  gradient methods for high-dimensional statistical recovery, in: Advances in
  Neural Information Processing Systems 23: 24th Annual Conference on Neural
  Information Processing Systems, NIPS, 2010, pp. 37--45.

\bibitem{huang2018constructive}
J.~Huang, Y.~Jiao, Y.~Liu, X.~Lu, A constructive approach to $\ell_0$ penalized
  regression, J. Mach. Learn. Res. 19~(10) (2018) 1--37.

\bibitem{garg2009gradient}
R.~Garg, R.~Khandekar, Gradient descent with sparsification: an iterative
  algorithm for sparse recovery with restricted isometry property, in:
  Proceedings of the 26th Annual International Conference on Machine Learning,
  ACM, 2009, pp. 337--344.

\end{thebibliography}
\end{document}